# ROBUST INFERENCE FOR UNIVARIATE PROPORTIONAL HAZARDS FRAILTY REGRESSION MODELS

By Michael R. Kosorok[1], Bee Leng Lee[2] and Jason P. Fine

*University of Wisconsin, Madison, National University of Singapore and University of Wisconsin, Madison*

We consider a class of semiparametric regression models which are one-parameter extensions of the Cox [*J. Roy. Statist. Soc. Ser. B* **34** (1972) 187–220] model for right-censored univariate failure times. These models assume that the hazard given the covariates and a random frailty unique to each individual has the proportional hazards form multiplied by the frailty. The frailty is assumed to have mean 1 within a known one-parameter family of distributions. Inference is based on a nonparametric likelihood. The behavior of the likelihood maximizer is studied under general conditions where the fitted model may be misspecified. The joint estimator of the regression and frailty parameters as well as the baseline hazard is shown to be uniformly consistent for the pseudo-value maximizing the asymptotic limit of the likelihood. Appropriately standardized, the estimator converges weakly to a Gaussian process. When the model is correctly specified, the procedure is semiparametric efficient, achieving the semiparametric information bound for all parameter components. It is also proved that the bootstrap gives valid inferences for all parameters, even under misspecification. We demonstrate analytically the importance of the robust inference in several examples. In a randomized clinical trial, a valid test of the treatment effect is possible when other prognostic factors and the frailty distribution are both misspecified. Under certain conditions on the covariates, the ratios of the regression parameters are still identifiable. The practical utility of the procedure is illustrated on a non-Hodgkin's lymphoma dataset.

Received January 2002; revised July 2003.

[1]Supported in part by Grant CA75142 from the National Cancer Institute.
[2]Supported in part by Grant R-155-000-022-112 from the National University of Singapore.

**1. Introduction.** An objective of many medical studies is a predictive model for survival. The Cox (1972) model is popular for such analyses, because of its theoretical properties and availability in software packages. Unfortunately, in many practical settings the phenomenon under study is quite complicated and the assumed model is at best a working approximation to the truth. Consider the Non-Hodgkin's Lymphoma Prognostic Factors Project (1993) which analyzed data from a collection of cancer clinical trials. A system was developed to classify patients according to baseline characteristics. The scheme employs a proportional hazards model with five influential covariates. The ordinal and continuous predictors are dichotomized for clinical interpretation. There are also important risk factors which are omitted, such as treatment center. Diagnostics show that the model fits poorly [Gray (2000)]. Furthermore, the survival estimates are quite biased by the misspecification.

There are several alternatives to the Cox model which might improve the fit. These include additive hazards regression models [Aalen (1978, 1980) and Lin and Ying (1994)], accelerated failure time models [Tsiatis (1990) and Wei, Ying and Lin (1990)] and time-varying coefficient models [Sargent (1997)]. Additional models have been developed for covariate-dependent heteroscedasticity and other departures from proportionality [Bagdonavičius and Nikulin (1999) and Hsieh (2001)].

Frailty models are a comparatively parsimonious representation which generalize the Cox model in a natural way. The misspecified and omitted covariates are described by an unobservable random variable $\log(W)$ unique to the linear predictor of each patient. Let $T$ be the failure time and $Z = \{Z(t), t \geq 0\}$ a $d \times 1$ vector process of possibly time-dependent covariates. Denote $\lambda\{t; \tilde{Z}(t), W\}$ as the hazard function of $T$ conditionally on $\tilde{Z}(t) = \{Z(s), s \leq t\}$ and $W$. The *proportional hazards frailty regression model* is

$$(1.1) \qquad \lambda\{t; \tilde{Z}(t), W\} = a(t) \exp\{\log(W) + \beta' Z(t)\},$$

where $\beta$ is a $d \times 1$ regression parameter, $a(t)$ is an unspecified base hazard function and prime ($'$) denotes transpose. Taking $f(w; \gamma)$ to be the Lebesgue density of a continuous frailty $W$, where $\gamma$ is an unknown scalar, yields a rich class of semiparametric models. This class excludes models with positive probability of $W = 0$. Examples in the class include the inverse Gaussian frailty [Hougaard (1984)], the positive stable frailty [Hougaard (1986)], the log-normal frailty [McGilchrist and Aisbett (1991)], the power variance frailty [Aalen (1988)], the uniform frailty [Lee and Klein (1988)] and the threshold frailty [Lindley and Singpurwalla (1986)]. While the one-parameter extension (1.1) of the Cox model is unlikely to address all misspecification, it is a point of departure. The objective of this paper is to provide a rigorous foundation for inference within this class of models, adopting the point of view that any model is at best a working approximation to the truth.

It is popular to let $W$ have a gamma distribution with mean 1 and variance $\gamma$. With time-independent covariates, the model is equivalent to the odds-rate regression [Dabrowska and Doksum (1988)],

$$h(T) = -\beta^T Z + \varepsilon_\gamma, \quad (1.2)$$

where $h(t)$ is an unspecified strictly monotone increasing function, and $\exp(\varepsilon_\gamma)$ has a Pareto($\gamma$) distribution. Fixing $\gamma = 0$ gives proportional hazards, while $\gamma = 1$ gives proportional odds [Bennett (1983)]. If $f(w; \gamma)$ has zero variance, then (1.1) reduces to the Cox model and efficient estimation of $\beta$ is straightforward with the partial likelihood [Andersen and Gill (1982)]. Estimation for the special case (1.2) with $\gamma$ known has been studied extensively [Pettitt (1982, 1984), Cheng, Wei and Ying (1995, 1997), Murphy, Rossini and van der Vaart (1997), Fine, Ying and Wei (1998), Scharfstein, Tsiatis and Gilbert (1998), Shen (1998) and Slud and Vonta (2004)]. When the parameter in the frailty distribution is unknown, these methods are not applicable. Asymptotic theory for maximum likelihood estimation of model (1.1) with clusters of size greater than or equal to 2 and shared gamma $W$ having unknown $\gamma$ was derived by Parner (1998). See Nielsen, Gill, Andersen and Sørensen (1992) and Murphy (1994, 1995) for related work. A unified theory for estimation in model (1.1) with uncorrelated data and general frailty distribution is not available.

In this paper the focus is on independent observations. The data setup and frailty model assumptions are given in Section 2. Bagdonavičius and Nikulin (1999) suggested ad hoc estimators for the parameters, but their large-sample properties were not established rigorously. The large-sample results in Parner (1998) can be adapted to the univariate gamma frailty setting with a correctly specified model, but do not apply to other frailty models and cannot be used to address model misspecification.

In Section 3, a likelihood-based procedure for model (1.1) is formally proposed, and the existence of likelihood maximizers and of both score and information operators is examined without requiring the model to be correctly specified. Section 4 establishes uniform consistency and weak convergence of the parameter estimators under mild identifiability conditions which ensure the uniqueness of the implied parameter corresponding to the maximizer of the asymptotic limit of the likelihood with respect to the true model. We also study properties of the estimators in settings where the model is not identifiable, as occurs when $\beta = 0$ and the frailty variance and baseline hazard are confounded. To our knowledge, this is the first attempt at asymptotic theory for misspecified nonparametric maximum likelihood estimation (NPMLE) for semiparametric survival regression models. White's (1982) work on robust parametric likelihood estimation is not directly applicable due to the presence of nonparametric components in (1.1). The closest related work is on asymptotic theory for the misspecified Cox model based on partial likelihood [Struthers and Kalbfleisch (1986), Lin and Wei (1989) and Sasieni (1993)]. However, these results do not apply to estimation based on full nonparametric likelihood.

Because the parametric and nonparametric components in (1.1) are estimated simultaneously, inference is complicated. Parner (1998) showed that the variance of the NPMLE for the gamma frailty model with cluster sizes greater than or equal to 2 can be consistently estimated by inverting a discrete observed information matrix. However, computing the required second derivatives can be difficult when the likelihood does not have a closed form, for example, with log-normal frailties. Furthermore, the limiting covariance function is extremely complicated and does not permit the construction of analytic confidence bands for functionals of the baseline hazard such as covariate-specific survival functions. The procedure we employ is to maximize the profile likelihood using a simple fixed-point algorithm for the baseline hazard motivated by the EM algorithm. We show that bootstrapping this procedure provides valid inference, including variance estimation and the construction of confidence bands for survival functions under model misspecification. The estimated survival probabilities may not be unbiased in large samples. However, it may be useful to interpret these quantities as minimizing the Kullback–Leibler discrepancy between the survival curves under the fitted and true models, conditionally on covariates.

In Section 5, the identifiability conditions given in Section 4 are shown to be satisfied when the model is correctly specified. We further verify that the estimators achieve the semiparametric variance bound [Sasieni (1992) and Bickel, Klaassen, Ritov and Wellner (1993), hereafter abbreviated BKRW] and are fully efficient for all model parameters. Section 6 evaluates the conditions of Section 4 under misspecification. In this case the estimators are still uniformly consistent and converge weakly, but inference must be based on an infinite-dimensional analogue to White's (1982) robust variance formula. Our contributions beyond Parner's (1998) work on the shared gamma frailty model are threefold. First, we study univariate data. Second, we allow general frailty distributions. Third, we permit misspecification.

The robust inferences are practically useful under some well-known misspecification mechanisms. To begin, we establish that when the true model has the form (1.1) but the choice of the distribution of $W$ is incorrect, the parameter estimate for a single covariate which is independent of one or more other misspecified covariates may be consistent up to sign. Note that all the covariates may be partially misspecified under mild restrictions. The setting applies in particular when assessing treatment effect in a randomized trial. Next we show that if the covariates $Z$ are correctly specified and $\mathrm{E}[b'Z|\beta_0'Z]$ is linear in $\beta_0'Z$ for all linear combinations $b'Z$, then the parameter estimates are consistent for $\alpha_2\beta_0$, where $\beta_0$ is the true regression parameter and $0 < \alpha_2 \in \mathbb{R}$. When the Cox model is used but the true model has frailty variance greater than 0, the estimated effect will be $\alpha_1 \beta_0$, where $\alpha_1 \in (0,1)$ and $\beta_0$ is the true effect. The conditional linearity assumption has been used by Li and Duan (1989) to establish similar robustness results under link function violations for parametric regression and for the Cox model based on partial likelihood but without censoring. Our results are applicable under independent censoring and are based on the full likelihood so that joint estimation of $\gamma$ and $A(t) = \int_0^t a(s)\,ds$ as well as $\beta$ is possible, which may be necessary for survival predictions.

While the focus of this paper is on independent survival times, many of the results and methods of proof are potentially applicable to multivariate failure time data. Furthermore, the fact that a proportional hazards model can be changed to a nonproportional hazards model by simply adding a frailty underscores the need to be careful when interpreting marginal inferences based on multivariate shared frailty models involving covariates. Nonproportional marginal hazards may not imply correlation of the failure times [Hougaard (2000)]. Extending the univariate results of this paper to the multivariate setting is an important topic for future research.

Several computational issues are discussed in Section 7, and the utility of the methods is illustrated on the lymphoma data in Section 8. All proofs are given in Section 9.

## 2. The data setup and frailty models.

2.1. *Data assumptions.* The data $\{X_i = (V_i, \delta_i, Z_i), i = 1, \ldots, n\}$ consist of $n$ i.i.d. realizations of $X = (V, \delta, Z)$, where $V = T \wedge C$, $\delta = \mathbf{1}\{T \leq C\}$, $x \wedge y$ denotes the minimum of $x$ and $y$, $\mathbf{1}\{B\}$ is the indicator of $B$ and $C$ is the right censoring time. The analysis is restricted to an interval $[0, \tau]$, where $\tau < \infty$. The covariate $Z \equiv \{Z(t), t \in [0, \tau]\}$ is assumed to be a caglad (left-continuous with right-hand limits) process with $Z(t) \in \mathbb{R}^d$, $t \in [0, \tau]$. We make the following additional assumptions:

(A1) $\mathrm{P}[C = 0] = 0$, $\mathrm{P}[C \geq \tau | Z] = \mathrm{P}[C = \tau | Z] > 0$ almost surely, and censoring is independent of $T$ given $Z$.
(A1') Condition (A1) is strengthened to require that $C$ and $Z$ are independent.
(A2) The true density of $T$ given $Z$, $f_0(t|Z)$, exists and is bounded over $t \in [0, \tau]$ almost surely, and $\mathrm{P}[T > \tau | Z] > 0$ almost surely.
(A3) The total variation of $Z(\cdot)$ on $[0, \tau]$ is $\leq m_0 < \infty$ almost surely, and $\mathrm{var}[Z(0+)]$ is positive definite, where for a real function $F$ with right-hand limits we define $F(t+) = \lim_{s \downarrow t} F(s)$.
(A3') Condition (A3) is strengthened to require that $Z = (Z_1, Z_2)$, where $Z_1 \in \mathbb{R}$ is time independent and $Z_1$ and $Z_2$ are stochastically independent.

(A3″) Condition A3 is strengthened to require that $Z$ is time independent and that $\mathrm{E}[b'Z|c'Z]$ is linear in $c'Z$ for all $b, c \in \mathbb{R}^d$.

Conditions (A1) and (A2) are somewhat standard for right-censored regression models, while condition (A3) is needed for both asymptotic normality in Section 4.3 and for parameter identifiability when the model is correctly specified in Section 5.1. The condition on $\mathrm{var}[Z(0+)]$ is similar to Parner's (1998) condition 2(g). The more restrictive assumptions (A1′), (A3′) and (A3″) are only used in Sections 6.2 and 6.3 for establishing robustness results under misspecification. An important example of when condition (A3′) holds is when $Z_1$ indicates treatment and treatment assignment has been randomized to ensure that $Z_2$, corresponding to other prognostic factors, is independent of $Z_1$. An important example of when condition (A3″) holds is when $Z$ is multivariate normal.

2.2. *Frailty model assumptions.* The frailty models we consider in this paper posit that the hazard function has the form (1.1). After integrating over $W$, the corresponding survival function at time $t$ given $Z = z$ becomes

$$
\begin{aligned}
S(t|z) &\equiv \mathrm{P}[T > t|Z = z] \\
&= \mathrm{E}\bigg[\exp\bigg\{-W \int_0^t e^{\beta' z(s)} \, dA(s)\bigg\}\bigg|Z = z\bigg] \\
&= \Lambda_\gamma\bigg\{\int_0^t e^{\beta' z(s)} \, dA(s)\bigg\},
\end{aligned}
\tag{2.1}
$$

where $W$ is continuous and independent of $Z$, $A(s) \equiv \int_0^s a(u)\,du$, $\Lambda_\gamma(t) \equiv \int_0^\infty e^{-wt} f(w;\gamma)\,dw$ is the Laplace transform of $W$ and $\gamma \in \mathbb{R}$ is an unknown parameter. At this point we are not assuming that the posited model agrees with the true conditional density $f_0$. The remainder of this section contains technical conditions on the frailty models, conditions (B), (C), (D1)–(D3), (E1) and (E2), which some readers may wish to skip over on the first reading.

We assume that the posited model consists of a family of frailty transforms $\{\Lambda_\gamma\}$ and a collection of indices $\{\psi = (\gamma, \beta, A)\}$, which satisfy:

(B) $\beta \in \overline{B}_0$, where $B_0 \subset \mathbb{R}^d$ contains 0 and is open, convex and bounded and where $\overline{B}$ denotes closure of a set $B$.
(C) There exist a constant $c_0$ and a continuous, decreasing function $\varepsilon_0 : [0,\infty) \mapsto (0, 3/4)$, so that $0 < c_0 < \varepsilon_0(0) < 3/4$, $\lim_{t\to\infty} \varepsilon_0(t) = 0$ and, for each positive $m, t < \infty$, there is an extension of $\Lambda_{(\cdot)}(\cdot) : [0,m) \times [0,t] \mapsto [0,1]$ having domain $[-\varepsilon_0(t), m] \times [0,t]$.

For the parametric component $\theta \equiv (\gamma, \beta)$, define the parameter set $\Theta \equiv (-c_0, m_1) \times B_0$ for some positive $m_1 < \infty$. In consequence of conditions (A3) and (B), let $1 \leq K_0 < \infty$ be the maximum possible value of 1 and both $e^{\beta'Z(t)}$ and $e^{-\beta'Z(t)}$ over $\beta \in \overline{B}_0$ and $t \in [0,\tau]$. Also let $\mathcal{A}$ be the collection of monotone increasing functions $A:[0,\tau] \mapsto [0,\infty)$, with $A(\tau) < \infty$, and define $\mathcal{A}_0$ to be the subset of $\mathcal{A}$ consisting of absolutely continuous functions with derivative satisfying $0 \leq a(t) < \infty$ for all $t \in [0,\tau]$. For $\gamma \in [-c_0, m_1]$, define $\mathcal{A}_{(\gamma)}$ so that $\mathcal{A}_{(\gamma)} = \mathcal{A}$ when $\gamma \geq 0$, and $\mathcal{A}_{(\gamma)} = \{A \in \mathcal{A}: A(\tau) \leq \varepsilon_0^{-1}(\gamma)/K_0\}$ when $\gamma < 0$. Also define $\mathcal{A}_{(\gamma)}^\circ = \mathcal{A}$ when $\gamma \geq 0$ and $\mathcal{A}_{(\gamma)}^\circ = \{A \in \mathcal{A}: A(\tau) < \varepsilon_0^{-1}(\gamma)/K_0\}$ when $\gamma < 0$. We can now define the index sets $\Psi \equiv \{\psi:(\gamma,\beta) \in \Theta, A \in \mathcal{A}_{(\gamma)}\} = \{\psi: A \in \mathcal{A}, \beta \in \overline{B}_0, \gamma \in [-\varepsilon(K_0 A(\tau)), m_1]\}$ and $\Psi_0 \equiv \{\psi:(\gamma,\beta) \in \Theta,\ A \in \mathcal{A}_0 \cap \mathcal{A}_{(\gamma)}^\circ\}$.

We need the following additional conditions on $\Lambda_\gamma$, where we define $\dot{\Lambda}_\gamma \equiv \partial \Lambda_\gamma(t)/(\partial t)$, $\ddot{\Lambda}_\gamma \equiv \partial \dot{\Lambda}_\gamma(t)/(\partial t)$, $G_\gamma \equiv -\log \Lambda_\gamma$, $\dot{G}_\gamma \equiv \partial G_\gamma(t)/(\partial t)$, $\ddot{G}_\gamma \equiv \partial \dot{G}_\gamma(t)/(\partial t)$, $\dddot{G}_\gamma \equiv \partial \ddot{G}_\gamma(t)/(\partial t)$, $G_\gamma^{(1)} \equiv \partial G_\gamma/(\partial \gamma)$, $\dot{G}_\gamma^{(1)} \equiv \partial G_\gamma^{(1)}(t)/(\partial t)$ and $\ddot{G}_\gamma^{(1)} \equiv \partial \dot{G}_\gamma^{(1)}(t)/(\partial t)$:

(D1) For each positive $t < \infty$, we have the following for all $\gamma \in [-\varepsilon_0(t), m_1]$ and $u \in [0,t]$: $\Lambda_\gamma(0+) = 1$, $\Lambda_\gamma(u) > 0$, $0 < -\dot{\Lambda}_\gamma(u) < \infty$ and $0 \leq \ddot{\Lambda}_\gamma(u) < \infty$; $\partial \dddot{G}_\gamma(u)/(\partial u)$, $\partial \ddot{G}_\gamma^{(1)}(u)/(\partial u)$, $\partial \ddot{G}_\gamma^{(1)}(u)/(\partial \gamma)$ and $\partial^2 \dot{G}_\gamma^{(1)}(u)/(\partial \gamma)^2$ exist and are bounded; $\dot{G}_\gamma(0+) = 1$, $\dot{G}_\gamma^{(1)}(0+) = 0$, $\ddot{G}_\gamma(u) \leq 0$ and $\ddot{G}_\gamma^{(1)}(0+) < 0$.

(D2) There exists a $c_1:(0,m_1] \mapsto (0,\infty]$ such that, for any sequence $\{\gamma_k\} \in [0,m_1]$ with $\gamma_k \to \gamma > 0$, $\limsup_{k\to\infty} \sup_{u \geq 0} u^{c_1(\gamma)} \Lambda_{\gamma_k}(u) < \infty$ and $\limsup_{k\to\infty} \sup_{u \geq 0} |u^{1+c_1(\gamma)} \dot{\Lambda}_{\gamma_k}(u)| < \infty$.

(D3) There exists a $c_2:[0,m_1] \mapsto (0,\infty]$, with $c_2(0) = \infty$, such that for all sequences $t_k \to \infty$ and $\{\gamma_k\} \in [-\varepsilon_0(t_k), m_1]$ with $\gamma_k \to \gamma \geq 0$, $\liminf_{k\to\infty} \inf_{u \in [0,t_k]} t_k \dot{G}_{\gamma_k}(u) \geq c_2(\gamma)$.

(E1) For all $\gamma \in [0, m_1]$ and all $t \in [0,\infty)$, $\dot{G}_\gamma(t) + t\ddot{G}_\gamma(t) > 0$ and

$$\frac{\ddot{G}_\gamma(t)}{\dot{G}_\gamma(t)} + t\left[\frac{\dddot{G}_\gamma(t)}{\dot{G}_\gamma(t)} - \left\{\frac{\ddot{G}_\gamma(t)}{\dot{G}_\gamma(t)}\right\}^2\right] \leq 0.$$

(E2) $\lim_{\gamma \downarrow 0} \mathrm{E}[(W-1)^2]/\gamma = 1$ and $\lim_{\gamma \downarrow 0} \mathrm{E}[|W-1|^3]/\gamma = 0$, where $W$ is a random variable with Laplace transform $\Lambda_\gamma$.

Conditions (D1)–(D3) are needed for uniform consistency and weak convergence of the estimators. Condition (D1) is also used for identifiability when the model is correctly specified. Conditions (E1) and (E2) are needed for identifiability under misspecification.

REMARK 1. Condition (C) ensures that $\gamma_0 = 0$ is an interior point. Parts of (D1)–(D3) are conditions on the moments of $W$. For (D1) this follows since $\dot{G}_\gamma(0+) = \mathrm{E}[W]$ and $\ddot{G}_\gamma(0+) = -\mathrm{var}[W]$. Condition (D2) is satisfied if there exists a continuous function $c_1:(0,m_1] \mapsto (0,\infty)$ such that $\mathrm{E}[W^{-c_1(\gamma)}] < \infty$ for all $\gamma \in (0,m_1]$.

2.3. *Examples of frailty models.* The following are instances of frailty transforms:

1. The *gamma frailty* has $\Lambda_\gamma(t) = (1+\gamma t)^{-1/\gamma}$.
2. The *inverse Gaussian frailty* [Hougaard (1984)] has $\Lambda_\gamma(t) = \exp\{-\gamma^{-1}[(1+2\gamma t)^{1/2} - 1]\}$.
3. The *log-normal frailty* [McGilchrist and Aisbett (1991)] has
$$\Lambda_\gamma(t) = \int_\mathbb{R} \exp\{-te^{\gamma^{1/2}v - \gamma/2}\}\phi(v)\,dv.$$
4. The *positive stable frailty* [Hougaard (1986)] has $\Lambda_\gamma(t) = \exp\{-t^\gamma\}$.
5. The IGG($\alpha$) *family of frailty transforms* has the form
$$\Lambda_\gamma = \exp\left\{-\frac{1-\alpha}{\alpha\gamma}\left[\left(1+\frac{\gamma t}{1-\alpha}\right)^\alpha - 1\right]\right\},$$
where $\alpha \in [0,1)$ is assumed known. The IGG(0) family is obtained by taking the limit as $\alpha \downarrow 0$.

REMARK 2. The IGG($\alpha$) family includes both the gamma frailty ($\alpha=0$) and the inverse Gaussian frailty ($\alpha=1/2$). The IG is for "inverse Gaussian" and the second G for "gamma."

REMARK 3. $\Lambda_\gamma$ and the functionals of $\Lambda_\gamma$ introduced above are defined at $\gamma=0$ by continuity. In all the above frailties, excepting the positive stable frailty, $\Lambda_0(t) = \lim_{\gamma\to 0}\Lambda_\gamma(t) = e^{-t}$, corresponding to the Cox model.

The following states that most of the stated frailty conditions are valid for a number of standard frailty families:

PROPOSITION 1. *Conditions* (C), (D1)–(D3) *and* (E2) *are satisfied by the gamma, inverse Gaussian, log-normal and* IGG($\alpha$), *for any fixed* $\alpha \in [0,1)$, *frailty distributions.*

REMARK 4. Verification of these conditions for the log-normal is hard technically since $\Lambda_\gamma$ does not have a closed form.

REMARK 5. Condition (E1) is easily verified for the gamma, inverse Gaussian and IGG($\alpha$) frailties, and has been validated numerically for the log-normal frailty for $\gamma \in [0, 4.62]$, corresponding to a frailty variance of 100. We conjecture that (E1) holds for the log-normal frailty for all $0 \leq \gamma < \infty$.

REMARK 6. For the positive stable frailty, conditions (D1), (D3), (E1) and (E2) are not satisfied but conditions (C) and (D2) are. For example, $\dot{G}_\gamma(0+) = \infty$ when $\gamma < 1$. Note also for this frailty that, when $Z$ is time independent, $-\log S(t|z) = e^{\gamma \beta' z} + A^\gamma(t)$, and the model is thus not identifiable.

### 3. Nonparametric maximum likelihood estimation.

3.1. *The estimator.* The nonparametric log-likelihood has the form

$$L_n(\psi) \equiv \mathbb{P}_n\bigg\{\int_0^\tau [\log \dot{G}_\gamma(H^\psi(s)) + \beta'Z(s) + \log a(s)]\,dN(s)$$
(3.1)
$$- G_\gamma(H^\psi(V))\bigg\},$$

where $N(t) \equiv \mathbf{1}\{V \leq t, \delta = 1\}$, $Y(t) \equiv \mathbf{1}\{V \geq t\}$, $H^\psi(t) \equiv \int_0^t e^{\beta'Z(s)}\,dA(s)$, $a \equiv dA/dt$ and $\mathbb{P}_n$ is the expectation with respect to the empirical probability measure. As discussed by Murphy, Rossini and van der Vaart (1997), the maximum likelihood estimator for $a$ does not exist, because any unrestricted maximizer of (3.1) puts mass only at observed failure times and is not a continuous hazard.

Instead, we compute the maximizer by profiling over $A$. This yields estimators for $\theta \equiv (\gamma, \beta)$ and $A$, but not $a$. The profile likelihood is $pL_n(\theta) \equiv \sup_{A \in \mathcal{A}_{(\gamma)}} L_n(\psi) = L_n(\theta, \hat{A}_\theta)$, where $\hat{A}_\theta \equiv \arg\max_{A \in \mathcal{A}_{(\gamma)}} L_n(\theta, A)$. Consider one-dimensional submodels for $A$, $t \mapsto A_t(\cdot) \equiv \int_0^{(\cdot)}\{1 + th(s)\}\,dA(s)$, where $(\cdot)$ denotes an argument ranging over $[0, \tau]$ and $h : [0, \tau] \mapsto \mathbb{R}$ is a bounded function. When the $\gamma$ component of $\theta$ is nonnegative, the upper limit of elements of $\mathcal{A}_{(\gamma)}$ is unconstrained. In this setting one may differentiate $L_n\{(\theta, A_t)\}$ with respect to $t$, where $h(s) = \mathbf{1}\{s \leq u\}$ and $u \in [0, \tau]$, and solve for $A$ with $t = 0$, since $A = A_{t=0}$. Hence, $\hat{A}_\theta$ solves

$$\hat{A}_\theta(u) = \int_0^u \bigg(\mathbb{P}_n\bigg[Y(s)e^{\beta'Z(s)}\bigg(\dot{G}_\gamma\{H^{\hat{\psi}_\theta}(V)\} - \delta\frac{\ddot{G}_\gamma\{H^{\hat{\psi}_\theta}(V)\}}{\dot{G}_\gamma\{H^{\hat{\psi}_\theta}(V)\}}\bigg)\bigg]\bigg)^{-1}$$
(3.2)
$$\times \mathbb{P}_n\{dN(s)\}$$
$$\equiv \int_0^u \{J_n^{\hat{\psi}_\theta}(s)\}^{-1}\mathbb{P}_n\{dN(s)\},$$

where $\hat{\psi}_\theta \equiv (\theta, \hat{A}_\theta)$.

Under model misspecification, it is possible that the best fit will occur for some $\gamma < 0$. In this case, $f(w; \gamma)$ will usually not be a density, even though the quantity $S(t|z)$ in (2.1) is a proper survival function provided $A(\tau)$ is not too large. Specifically, all $A \in \mathcal{A}_{(\gamma)}$ must satisfy $A(\tau) \leq \varepsilon_0^{-1}(\gamma)/K_0$. Under this constraint, one may differentiate $L_n\{(\theta, A_t)\}$ with respect to $t$, where $h(s) = \mathbf{1}\{v < s \leq u\} - [A(u) - A(v)]/A(\tau)$ and $u, v \in [0, \tau]$, take $t = 0$, let $v \uparrow u$ and solve for $A$. This yields

(3.3) $$\hat{A}_\theta(u) = \int_0^u \{J_n^{\hat{\psi}_\theta}(s) + \rho_n(\hat{\psi}_\theta)\}^{-1}\mathbb{P}_n\{dN(s)\},$$

where

$$\rho_n(\psi) \equiv \frac{\mathbb{P}_n \delta - \int_0^\tau J_n^\psi(s)\, dA(s)}{A(\tau)}.$$

By considering one-dimensional submodels $\{A_s\}$ with $h(s) = -1$, the fact that the derivative of $L_n(\theta, A_s)$ at $s = 0$ is nonpositive implies that $\rho_n(\hat{\psi}_\theta) \geq 0$. Thus, for all $\theta \in \overline{\Theta}$, $\hat{A}_\theta$ has the form given in (3.3), with $\rho_n(\hat{\psi}_\theta) > 0$ only when $\gamma < 0$ and $\hat{A}_\theta(\tau) = \varepsilon_0^{-1}(\gamma)/K_0$, and $\rho_n(\hat{\psi}_\theta) = 0$ otherwise.

The same maximizer occurs with $\Delta A$ in place of $a$ in $L_n(\psi)$, where $\Delta A(s) \equiv A(s) - A(s-)$ and $A(s-) \equiv \lim_{t \uparrow s} A(t)$. That is, one maximizes $L_n(\psi)$ over all $A$ with jumps at the observed failure times. We denote by $\tilde{L}_n(\psi)$ the log-likelihood expression with $\Delta A$ in place of $a$. The nonparametric maximum likelihood estimator (NPMLE) is $\hat{\psi}_n \equiv (\hat{\theta}_n, \hat{A}_{\hat{\theta}_n})$, where $\hat{\theta}_n \equiv \arg\max_{\theta \in \overline{\Theta}} pL_n(\theta)$. Equivalently, $\hat{\psi}_n = \arg\max_{\psi \in \Psi} \tilde{L}_n(\psi)$.

We have the following existence result.

PROPOSITION 2. *Under conditions* (A1)–(A3), (B), (C) *and* (D1)–(D3), *and provided* $\max_{1 \leq i \leq n} \delta_i > 0$, *then for some* $1 < M < \infty$ *and all* $\theta \in \overline{\Theta}$ *there exist maximizers* $\hat{A}_\theta$, *and all such maximizers satisfy* (3.3) *and* $1/M \leq \hat{A}_\theta \leq M$.

REMARK 7. Proposition 2 implies the existence of an NPMLE $\hat{\psi}_n$ as a consequence of the compactness of $\overline{\Theta}$. However, the proposition says nothing about uniqueness of the NPMLE.

For a limiting value of the NPMLE to exist, it is necessary (but not sufficient) that the $M$ in Proposition 2 does not go to $\infty$ as $n \to \infty$. However, a significantly stronger result can be obtained. Define $\mathcal{K}_M = \{A \in \mathcal{A}: 1/M \leq A(\tau) \leq M, \sup_{t \in [0,\tau]} a(t) \leq M\}$ and, for each $\varepsilon > 0$, define $\mathcal{K}_M^\varepsilon = \{A \in \mathcal{A}: \sup_{t \in [0,\tau]} |A(t) - \tilde{A}(t)| \leq \varepsilon \text{ for some } \tilde{A} \in \mathcal{K}_M\}$. Note that $\mathcal{K}_M$ is compact for each $1 < M < \infty$. Let $P_*$ denote inner probability. We have the following result.

THEOREM 1. *Assume conditions* (A1)–(A3), (B), (C) *and* (D1)–(D3). *Then, for each* $\eta > 0$, *there exist some* $1 < M < \infty$ *such that* $\lim_{\varepsilon \downarrow 0} P_*(\{\hat{A}_\theta : \theta \in \overline{\Theta}\} \in \mathcal{K}_M^\varepsilon \ \forall n \text{ large enough}) > (1 - \eta)$.

REMARK 8. Theorem 1 implies that all sequences of NPMLE's have convergent subsequences and that the resulting limit points for $\hat{A}_{\hat{\theta}_n}$ have bounded derivatives almost surely. Consistency will then follow from identifiability of the model. Moreover, when only some of the parameters are identifiable, consistency of the identifiable parameters will also follow. The important example of estimation of $\beta$ when the survival distribution does not depend on covariates is discussed in Section 4.2.

3.2. *Kullback–Leibler information.* We now establish properties of the Kullback–Leibler information. Let $p_\psi(v,e|z) \equiv f^e(v|z) S^{1-e}(v|z)$, where $f(t|z) \equiv -\partial S(t|z)/(\partial t)$, $S(t|z) = \exp\{-G_\gamma(H^\psi(t))\}$ as defined in (2.1). For each $\theta \in \overline{\Theta}$, let $A_\theta \equiv \arg\max_{A \in \mathcal{A}_{(\gamma)} \cap \mathcal{A}_0} P_0 \log(p_\psi)$ and $\psi_\theta \equiv (\theta, A_\theta)$. We have the following result.

LEMMA 1. *Under conditions* (A1)–(A3), (B), (C) *and* (D1)–(D3), *and for some* $1 < M < \infty$, $A_\theta \in \mathcal{K}_M$ *for all* $\theta \in \overline{\Theta}$.

REMARK 9. Lemma 1 tells us that, even without any model identifiability, all possible Kullback–Leibler maximizers lie in a compact set. Questions about consistency can thus be reduced to questions about identifiability or partial identifiability as mentioned in Remark 8.

For each $\theta \in \overline{\Theta}$, let $\psi_\theta = (\theta, A_\theta)$ and $\hat{\psi}_\theta = (\theta, \hat{A}_\theta)$. For any $\psi_1, \psi_2$ in the subset of $\Psi$ where $A_1$ and $A_2$ have jumps only at observed failure times, define the empirical Kullback–Leibler information $I_n(\psi_1, \psi_2) \equiv \tilde{L}_n(\psi_1) - \tilde{L}_n(\psi_2)$, and, for any $\psi_1, \psi_2$ in the subset of $\Psi$ for which the derivatives $a_1$ and $a_2$ exist, define the Kullback–Leibler information $I_0(\psi_1, \psi_2) \equiv P_0 \log(p_{\psi_1}/p_{\psi_2})$. The following theorem establishes, in the profile context, an important asymptotic equivalence between $I_n$ and $I_0$.

THEOREM 2. *Under conditions* (A1)–(A3), (B), (C) *and* (D1)–(D3),
$$\sup_{\theta_1, \theta_2 \in \overline{\Theta}} |I_n(\hat{\psi}_{\theta_1}, \hat{\psi}_{\theta_2}) - I_0(\psi_{\theta_1}, \psi_{\theta_2})| \to 0$$
*outer almost surely, as* $n \to \infty$.

REMARK 10. While Proposition 2 and Lemma 1 establish the existence of the profile maximizers $\hat{A}_\theta$ and $A_\theta$, uniqueness is not established. However, Theorem 2 tells us that all members of the equivalence class $\hat{A}_\theta$ are asymptotically equivalent to all members of the equivalence class $A_\theta$ in terms of Kullback–Leibler information. Thus, model identifiability immediately implies asymptotic uniqueness.

3.3. *Score and information operators.* In this section we derive the score and information operators. These play a key role in the weak convergence results presented in later sections. For each $\psi \in \Psi$ with $A$ having bounded derivative, define the one-dimensional submodels $t \mapsto \psi_t \equiv \psi + t\{h_1, h_2, \int_0^{(\cdot)} h_3(s) \, dA(s)\}$, where $(h_1, h_2, h_3) \in H_r$ for some $r < \infty$ and where $H_r$ is the space of elements $h = (h_1, h_2, h_3)$ such that $h_1 \in \mathbb{R}$, $h_2 \in \mathbb{R}^d$, $h_3$ is a cadlag (right-continuous with left-hand limits) function and $|h_1| + \sqrt{h_2' h_2} + \|h_3\|_v \leq r$, with $\|\cdot\|_v$ being the total variation norm. Let $H_\infty = \bigcup_{0 < r < \infty} H_r$. Since $\psi$ can be represented as a functional on $H_r$ of the form $\psi(h) = h_1 \gamma + h_2' \beta + \int_0^\tau h_3(s) \, dA(s)$, the parameter space $\Psi$ is then a subset of $\ell^\infty(H_r)$ with norm $\|\psi\|_{(r)} \equiv \sup_{h \in H_r} |\psi(h)|$, where $\ell^\infty(B)$ is the space of bounded functionals on $B$. For $\psi \in \Psi$ and $g, h \in H_r$, define

$$\psi^g(h) \equiv g_1 h_1 + g_2' h_2 + \int_0^\tau g_3(s) h_3(s) \, dA(s).$$

Note that $H_1$ is rich enough to extract all components of $\psi$ since $H_1$ includes $\{h : h_1 = 1, h_2 = h_3 = 0\} \cup \{h : h_2 = 1, h_1 = h_3 = 0\} \cup \{h : h_1 = h_2 = 0, h_3(s) = \mathbf{1}\{s \leq t\}, t \in [0, \tau]\}$. Let

$$\begin{aligned}
U_n^\tau(\psi)(h) &\equiv \frac{\partial}{\partial t} L_n(\psi_t) \Big|_{t=0} \\
&= \mathbb{P}_n \Bigg\{ \bigg[ \delta \frac{\dot{G}_\gamma^{(1)}(H^\psi(V \wedge \tau))}{\dot{G}_\gamma(H^\psi(V \wedge \tau))} - G_\gamma^{(1)}(H^\psi(V \wedge \tau)) \bigg] h_1 \\
&\quad + \int_0^\tau [Z'(s) h_2 + h_3(s)] \, dN(s) \\
&\quad + \bigg[ \delta \frac{\ddot{G}_\gamma(H^\psi(V \wedge \tau))}{\dot{G}_\gamma(H^\psi(V \wedge \tau))} - \dot{G}_\gamma(H^\psi(V \wedge \tau)) \bigg] \\
&\quad \times \int_0^\tau Y(s) e^{\beta' Z(s)} [Z'(s) h_2 + h_3(s)] \, dA(s) \Bigg\} \\
&\equiv \mathbb{P}_n U^\tau(\psi)(h).
\end{aligned}$$

(3.4)

This score operator can easily be extended so that the bounded derivative restriction on $A$ is unnecessary. The operator has expectation $U_0^\tau(\psi) \equiv P_0 U^\tau(\psi)$. The dependence on $\tau$ will be needed later.

The Gâteaux derivative of $U_0^\tau(\psi)(h)$ at $\psi_1 \in \Psi$ exists and is obtained by differentiating the score operator for the submodels $t \mapsto \psi_1 + t\psi$. This derivative is

$$-\dot{U}_{\psi_1}(\psi)(h) \equiv -\frac{\partial}{\partial t} U_0^\tau(\psi_1 + t\psi)(h) \Big|_{t=0} = \psi(\sigma_{\psi_1}(h)),$$

where the operator $\sigma_\psi : H_\infty \mapsto H_\infty$ is

$$\sigma_\psi(h) = \begin{pmatrix} \sigma_\psi^{11} & \sigma_\psi^{12} & \sigma_\psi^{13} \\ \sigma_\psi^{21} & \sigma_\psi^{22} & \sigma_\psi^{23} \\ \sigma_\psi^{31} & \sigma_\psi^{32} & \sigma_\psi^{33} \end{pmatrix} \begin{pmatrix} h_1 \\ h_2 \\ h_3 \end{pmatrix}$$

(3.5)
$$\equiv P_0 \begin{pmatrix} \hat{\sigma}_\psi^{11} & \hat{\sigma}_\psi^{12} & \hat{\sigma}_\psi^{13} \\ \hat{\sigma}_\psi^{21} & \hat{\sigma}_\psi^{22} & \hat{\sigma}_\psi^{23} \\ \hat{\sigma}_\psi^{31} & \hat{\sigma}_\psi^{32} & \hat{\sigma}_\psi^{33} \end{pmatrix} \begin{pmatrix} h_1 \\ h_2 \\ h_3 \end{pmatrix}$$

$$= P_0 \hat{\sigma}_\psi(h).$$

The operators $\sigma_\psi^{jk} = P_0 \hat{\sigma}_\psi^{jk}$, for $1 \leq j, k \leq 3$, are well defined and bounded, where

$$\hat{\sigma}_\psi^{11}(h_1) = \hat{\xi}_\psi^{(1)} h_1,$$

$$\hat{\sigma}_\psi^{21}(h_1) = \hat{\xi}_\psi^{(2)} \int_0^\tau Z(s) Y(s) e^{\beta' Z(s)} \, dA(s) h_1,$$

$$\hat{\sigma}_\psi^{31}(h_1)(t) = \hat{\xi}_\psi^{(2)} Y(t) e^{\beta' Z(t)} h_1,$$

$$\hat{\sigma}_\psi^{12}(h_2) = \hat{\xi}_\psi^{(2)} \int_0^\tau Z'(s) h_2 Y(s) e^{\beta' Z(s)} \, dA(s),$$

$$\hat{\sigma}_\psi^{13}(h_3) = \hat{\xi}_\psi^{(2)} \int_0^\tau h_3(s) Y(s) e^{\beta' Z(s)} \, dA(s),$$

$$\hat{\sigma}_\psi^{22}(h_2) = \hat{\xi}_\psi^{(0)} \int_0^\tau Z(s) Z'(s) h_2 Y(s) e^{\beta' Z(s)} \, dA(s)$$
$$+ \hat{\xi}_\psi^{(3)} \int_0^\tau Z(s) Y(s) e^{\beta' Z(s)} \, dA(s) \int_0^\tau Z'(s) h_2 Y(s) e^{\beta' Z(s)} \, dA(s),$$

$$\hat{\sigma}_\psi^{23}(h_3) = \hat{\xi}_\psi^{(0)} \int_0^\tau Z(s) h_3(s) Y(s) e^{\beta' Z(s)} \, dA(s)$$
$$+ \hat{\xi}_\psi^{(3)} \int_0^\tau Z(s) Y(s) e^{\beta' Z(s)} \, dA(s) \int_0^\tau h_3(s) Y(s) e^{\beta' Z(s)} \, dA(s),$$

$$\hat{\sigma}_\psi^{32}(h_2)(t) = \hat{\xi}_\psi^{(0)} Z'(t) h_2 Y(t) e^{\beta' Z(t)}$$
$$+ \hat{\xi}_\psi^{(3)} Y(t) e^{\beta' Z(t)} \int_0^\tau Z'(s) h_2 Y(s) e^{\beta' Z(s)} \, dA(s),$$

$$\hat{\sigma}_\psi^{33}(h_3)(t) = \hat{\xi}_\psi^{(0)} h_3(t) Y(t) e^{\beta' Z(t)}$$
$$+ \hat{\xi}_\psi^{(3)} Y(t) e^{\beta' Z(t)} \int_0^\tau h_3(s) Y(s) e^{\beta' Z(s)} \, dA(s)$$

and where

$$\hat{\xi}_\psi^{(0)} = \dot{G}_\gamma(H^\psi(V \wedge \tau)) - \delta \frac{\ddot{G}_\gamma(H^\psi(V \wedge \tau))}{\dot{G}_\gamma(H^\psi(V \wedge \tau))},$$

$$\hat{\xi}_\psi^{(1)} = G_\gamma^{(2)}(H^\psi(V \wedge \tau))$$
$$- \delta \left[ \frac{\dot{G}_\gamma^{(2)}(H^\psi(V \wedge \tau))}{\dot{G}_\gamma(H^\psi(V \wedge \tau))} - \left\{ \frac{\dot{G}_\gamma^{(1)}(H^\psi(V \wedge \tau))}{\dot{G}_\gamma(H^\psi(V \wedge \tau))} \right\}^2 \right],$$

$$\hat{\xi}_\psi^{(2)} = \dot{G}_\gamma^{(1)}(H^\psi(V \wedge \tau))$$
$$- \delta \left[ \frac{\ddot{G}_\gamma^{(1)}(H^\psi(V \wedge \tau))}{\dot{G}_\gamma(H^\psi(V \wedge \tau))} - \frac{\dot{G}_\gamma^{(1)}(H^\psi(V \wedge \tau))\ddot{G}_\gamma(H^\psi(V \wedge \tau))}{\{\dot{G}_\gamma(H^\psi(V \wedge \tau))\}^2} \right]$$

and

$$\hat{\xi}_\psi^{(3)} = \ddot{G}_\gamma(H^\psi(V \wedge \tau)) - \delta \left[ \frac{\dddot{G}_\gamma(H^\psi(V \wedge \tau))}{\dot{G}_\gamma(H^\psi(V \wedge \tau))} - \left\{ \frac{\ddot{G}_\gamma(H^\psi(V \wedge \tau))}{\dot{G}_\gamma(H^\psi(V \wedge \tau))} \right\}^2 \right],$$

where we also define $\dot{G}_\gamma^{(2)} \equiv \partial \dot{G}_\gamma^{(1)}(t)/(\partial \gamma)$.

To use the $Z$-estimator master theorem to obtain weak convergence in Section 4.3, the Gâteaux differentiability of $U_0^\tau$ needs to be strengthened to Fréchet differentiability. Accordingly, we have the following result.

LEMMA 2. *Under conditions* (A1)–(A3), (B), (C) *and* (D1)–(D3), *and for any $\psi_1 \in \Psi$, the operator $\psi \mapsto U_0^\tau(\psi)$ is Fréchet differentiable at $\psi_1$, with derivative $\psi(\sigma_{\psi_1}(h))$.*

## 4. General results.

4.1. *Additional assumptions.* Let $S_0(t|z) = \int_t^\infty f_0(s|z)\,ds$, where $f_0$ is as defined in condition (A2). Denote $p_0(v,e|z) = f_0^e(v|z)S_0^{(1-e)}(v|z)$ and let $\nu(v,e,z)$ be the implicitly defined measure for $(V,\delta,Z)$ such that the true expectation of $g(X)$, denoted $P_0 g$, can be written as $\int_\mathcal{X} gp_0 \,d\nu$, where $\mathcal{X}$ is the sample space for $X$ and $g$ is measurable. Recall that the operator $\sigma_\psi$ was defined in (3.5) for all $\psi \in \Psi$. We make the following assumptions about the relationship between the posited frailty model and the true distribution:

(F) $P_0 \log(p_\psi/p_0)$ has a unique maximum over $\psi \in \Psi_0$ at $\psi_* = (\gamma_*, \beta_*, A_*) \in \Psi_0$.
(G) $\sigma_{\psi_*}: H_\infty \mapsto H_\infty$ is one-to-one.

REMARK 11. Assumption (F) is analogous to assumption A3(b) of White (1982) and is required for consistency, while condition (G) is analogous to assumption A6(b) of White (1982) and is required for asymptotic normality. The lack of convexity of the Kullback–Leibler information in the posited frailty models generally prevents assumptions (F) and (G) from being direct consequences of the other conditions, except when the frailty model is correctly specified (Section 5) or when the true model is not too far from a member of the posited frailty model (Section 6.1).

REMARK 12. With a misspecified frailty model, the existence of the implicitly defined $\psi_*$ does not guarantee its meaningfulness. In general, $p_{\psi_*} \neq p_0$ and $\psi_* = \psi_0$ only when $p_{\psi_0} = p_0$. We show in Sections 6.2 and 6.3 that when $p_\psi$ is misspecified but assumption (F) holds, some of the components of $\psi_*$ may sometimes be useful for inference about $p_0$.

4.2. *Consistency.* The theorem we now present establishes the consistency of $\hat{\psi}_n$ under the identifiability assumed in (F).

THEOREM 3. *Under the conditions of Proposition* 2 *and condition* (F), $\hat{\psi}_n$ *converges outer almost surely to* $\psi_*$ *in the uniform norm.*

The following gives us the consistency of $\hat{\beta}_n$ under an important partial identifiability setting not requiring condition (F).

PROPOSITION 3. *Assume the conditions of Proposition* 2 *and that* $f_0(t|Z = z_1) = f_0(t|Z = z_2)$ *for all* $t \in [0, \tau]$ *and all possible values* $z_1$ *and* $z_2$ *of the covariate process. Then* $\hat{\beta}_n$ *converges outer almost surely to* 0.

REMARK 13. An innovation in the proofs of Theorem 3 and Proposition 3 is that the existence and asymptotic boundedness of $\hat{A}_n$ is established even when the model is misspecified or when condition (F) may not hold. This was shown in Proposition 2 and Theorem 1, where it was demonstrated that the asymptotic boundedness and equicontinuity of $\hat{A}_n$ depends only on the structure of the data, the posited model and the general condition (A2), but does not depend on any other aspects of the underlying true distribution.

4.3. *Asymptotic normality.* We use Hoffmann–Jørgensen weak convergence as described in van der Vaart and Wellner (1996) (hereafter abbreviated VW). We have the following result.

THEOREM 4. *Under the conditions of Theorem* 3 *and condition* (G), $\sqrt{n}(\hat{\psi}_n - \psi_*)$ *is asymptotically linear, with influence function* $\tilde{\ell}(h) = U^\tau(\psi_*) \times (\sigma_{\psi_*}^{-1}(h))$, $h \in H_1$, *converging weakly in the uniform norm to a tight, mean-zero Gaussian process* $\mathcal{Z}_*$ *with covariance* $V_*(g, h) = \mathrm{E}[\tilde{\ell}(g)\tilde{\ell}(h)]$, $g, h \in H_1$.

REMARK 14. In the proof in Section 9, the problem of establishing weak convergence can be cleanly divided into establishing properties of the data and fitted model (1.1), based on conditions (A1)–(A3), (B), (C) and (D1)–(D3), and establishing properties of the Kullback–Leibler discrepancy $P_0 \log(p_\psi/p_0)$, based on conditions (F) and (G), which involves the true distribution of the censoring and covariates.

4.4. *The bootstrap.* The usual nonparametric bootstrap resamples with replacement from the observed data. A disadvantage is that ties can arise with censored survival data. We propose an alternative weighted bootstrap. In each bootstrap sample one generates $n$ independent and identically distributed nonnegative weights $\zeta_1, \ldots, \zeta_n$, with mean and variance 1 and with $\int_0^\infty \sqrt{P[\zeta_1 > x]}\, dx < \infty$. Each weight is divided by the average weight (rejecting samples with all 0's) to obtain "standardized weights" $\zeta_1^\circ, \ldots, \zeta_n^\circ$ which sum to $n$. Distributions satisfying the moment conditions include the unit exponential and the Poisson with mean 1. For the nonparametric bootstrap the weights $\zeta_1^{\cdot}, \ldots, \zeta_n^{\cdot}$ are generated from a multinomial distribution with $E\zeta_i^{\cdot} = 1$, $i = 1, \ldots, n$, and $\sum_{i=1}^n \zeta_i^{\cdot} = n$.

For a known function $f$, let $\mathbb{P}_n^\circ f(V, \delta, Z; \psi) \equiv n^{-1} \sum_{i=1}^n \zeta_i^\circ f(V_i, \delta_i, Z_i; \psi)$ define the weighted empirical measure $\mathbb{P}_n^\circ$. The weighted bootstrap estimate $\hat{\psi}_n^\circ$ is computed by substituting $\mathbb{P}_n^\circ$ for $\mathbb{P}_n$ in the expressions in Section 2.2 and maximizing over $\psi$. Note that $\mathbb{P}_n^{\cdot}$ is defined similarly to $\mathbb{P}_n^\circ$ with the weights $\zeta_1^{\cdot}, \ldots, \zeta_n^{\cdot}$ in place of $\zeta_1^\circ, \ldots, \zeta_n^\circ$. The nonparametric bootstrap estimate $\hat{\psi}_n^{\cdot}$ is computed by using $\mathbb{P}_n^{\cdot}$ in place of $\mathbb{P}_n$ in Section 2.2.

The following result establishes the validity of both the nonparametric and the weighted bootstraps.

COROLLARY 1. *Under the assumptions of Theorem 4, the conditional bootstrap of $\hat{\psi}_n$, based either on $\hat{\psi}_n^{\cdot}$ or $\hat{\psi}_n^\circ$, is asymptotically consistent for the limiting process $\mathcal{Z}_*$. That is, $\sqrt{n}(\hat{\psi}_n^{\cdot} - \hat{\psi}_n)$ and $\sqrt{n}(\hat{\psi}_n^\circ - \hat{\psi}_n)$ are asymptotically measurable,*

(i) $\sup_{g \in BL_1} |E_{\cdot} g(\sqrt{n}(\hat{\psi}_n^{\cdot} - \hat{\psi}_n)) - Eg(\mathcal{Z}_*)| \to 0$ *in outer probability, and*
(ii) $\sup_{g \in BL_1} |E_\circ g(\sqrt{n}(\hat{\psi}_n^\circ - \hat{\psi}_n)) - Eg(\mathcal{Z}_*)| \to 0$ *in outer probability,*

*where $BL_1$ is the space of functions mapping $\mathbb{R}^{d+1} \times \ell^\infty([0, \tau]) \mapsto \mathbb{R}$ with Lipschitz norm $\leq 1$, and conditional on the data $E_{\cdot}$ and $E_\circ$ are expectations over the multinomial and standardized weights, respectively.*

REMARK 15. While the choice of $\{\zeta_i\}$ in the weighted bootstrap has no effect asymptotically, the rate of convergence may be affected. Newton and Raftery (1994) discuss different choices in the context of parametric maximum likelihood. They demonstrate that unit exponential weights, which are Dirichlet after standardizing, perform well. Our own experience is that exponential weights also work well for semiparametric inference. A detailed analysis of the distribution of the weights is beyond the scope of this paper.

REMARK 16. An advantage of using $Z$-estimator theory for establishing weak convergence of estimators for likelihood inference under possible model misspecification is that consistency of the bootstrap is essentially an immediate consequence of the influence function being $P_0$-Donsker.

**5. Results under correctly specified model.** The focus of this section is on the behavior of $\hat{\psi}_n$ when the frailty regression model is correctly specified. Accordingly, $\psi_* = \psi_0$ throughout this section. In Section 5.1 we establish that the identifiability condition [condition (F)] holds. In Section 5.2 the injectiveness of $\sigma_{\psi_0}$ [condition (G)] is established and shown to imply that $\hat{\psi}_n$ is both regular and efficient. In addition to assuming that $\Lambda_\gamma$ is correctly specified, we make the following assumption:

(H) $\gamma_0 \in [0, m_1)$, $0 \neq \beta_0 \in B_0$ and $A_0 \in \mathcal{A}_0$ with $a_0 > 0$.

REMARK 17. Condition (H) is also assumed by Parner (1998). When $\beta_0 = 0$, the survival function $S_0(t)$ does not depend on covariates, and we have the situation considered in Proposition 3. Moreover, we have that for any $\gamma \geq 0$, there exists an $A_\gamma \in \mathcal{A}_0$ so that $S_0(t) = \Lambda_\gamma(A_\gamma(t))$. Thus, $\gamma$ and $A$ are not identifiable when $\beta_0 = 0$.

5.1. *Identifiability.* Nonparametric identifiability of the mixed proportional hazards model, under the assumption that $\beta_0' Z$ takes on at least two distinct values, has been established for right-, left- and double-censored data with finite-mean frailties by Kortram, van Rooij, Lenstra and Ridder (1995). For earlier related work, see also Heckman and Taber (1994), Elbers and Ridder (1982) and Heckman and Singer (1984). In our case, $\Lambda_\gamma$ is parametric rather than completely unspecified and may not have the interpretation of being the Laplace transform of a frailty when $\gamma < 0$. The following proposition establishes uniqueness of the model.

PROPOSITION 4. *Under conditions* (A1)–(A3), (B), (C), (D1)–(D3) *and* (H), *model* (2.1) *is identifiable over* $\Psi_0$, *and thus condition* (F) *is satisfied.*

REMARK 18. The monotonicity in $\gamma$ of $\ddot{G}_\gamma(0+)$, where $\ddot{G}_\gamma(t) \equiv -\partial^2 \log \Lambda_\gamma(t)/(\partial t)^2$, as given in condition (D1), is the key to establishing identifiability of the extended Laplace transform. Since $-\ddot{G}_\gamma(0+)$ is the variance of $W$ when $\gamma \geq 0$, this is the same as requiring that $\text{var}[W]$ be a monotone function of $\gamma$. The positive stable frailty model violates this condition and is not identifiable without clustered data as noted in Remark 6 above. Because of Proposition 1, the gamma, inverse Gaussian, log-normal and IGG($\alpha$) frailties are identifiable.

5.2. *Efficiency.* The main result of this section is as follows.

THEOREM 5. *The information operator $\sigma_{\psi_0}$ is one-to-one. Thus, $\sigma_{\psi_0}$ is continuously invertible, condition* (G) *is satisfied and $\hat{\psi}_n$ is a regular and efficient estimator of $\psi_0$ when $\gamma_0 \geq 0$ and censoring is uninformative of $\psi$. The limiting covariance for $\sqrt{n}(\hat{\psi}_n - \psi_0)$ is $\psi_0^g(\sigma_{\psi_0}^{-1}(h))$, $g, h \in H_1$.*

REMARK 19. For the shared gamma frailty regression model, Parner (1998) suggests inference based on estimating the covariance. This is done through first estimating $\sigma_{\psi_0}$ by plugging in $\hat{\psi}_n$ for $\psi_0$ and then inverting, considering only the parameters $\gamma$, $\beta$ and $\Delta A$ at observed failure times. Because this approach may be difficult to implement with general $f(w; \gamma)$ and does not readily enable the construction of confidence bands, the bootstrap is recommended for inference. By Corollary 1, Theorem 5 implies that the bootstrap will yield valid inferences.

REMARK 20. The proof of Theorem 5 draws heavily on the tangent set $H_r$, as defined in the proof of Theorem 4. The issue is showing that, for $h \in H_r$, $\sigma_{\psi_0}(h) = 0$ implies $h = 0$. This gives that $\sigma_{\psi_0}$ is continuously invertible and onto, and thus the influence function is contained in the closed linear span of the score operator, yielding the given covariance. Regularity and efficiency then follow from Theorems 5.2.3 and 5.2.1 of BKRW.

REMARK 21. An alternative estimator to $\hat{\psi}_n$ is to take

$$\tilde{\psi}_n = \begin{cases} \hat{\psi}_n, & \text{if } \hat{\gamma}_n > 0, \\ (0, \tilde{\beta}_n, \tilde{A}_n), & \text{otherwise,} \end{cases}$$

where $\tilde{\beta}_n, \tilde{A}_n$ are the estimates based on the Cox model. It is not difficult to show that, when $\gamma_0 > 0$, $\tilde{\psi}_n$ has the same limiting behavior as $\hat{\psi}_n$, but when $\gamma_0 = 0$, the limiting distribution of $\sqrt{n}(\tilde{\psi}_n - \psi_0)$ is a mixture of the limiting distribution of $\sqrt{n}(\hat{\psi}_n - \psi_0)$ and the limiting distribution under the Cox model (with a 0 in the $\gamma$ component), each with probability 1/2. This alternative estimator is thus more precise when $\gamma_0 = 0$. It also follows without difficulty that the conditional limit law of the bootstrap which imitates this estimation procedure is equal to the limit law of $\sqrt{n}(\tilde{\psi}_n - \psi_0)$ in the sense of Corollary 1.

**6. Results under model misspecification.** In this section, we examine conditions under which the model is misspecified but the parameter estimates are consistent and asymptotically Gaussian, and some of the components of the estimated quantity may be interpreted via $p_0$. In Section 6.1 we demonstrate that if the posited conditional survival distribution—based on the chosen frailty transform $\Lambda_\gamma$—is not too badly misspecified, then conditions (F) and (G) are satisfied under certain restrictions on the index set $\Psi$. In Section 6.2, we examine the effect of testing for the effect of a single covariate with misspecification. In Section 6.3 we study $\gamma$ and $\beta$ under misspecification with structural requirements on the covariates.

6.1. *Existence of unique Kullback–Leibler maximizers.* Define

$$\Psi_M = \{\psi = (\gamma, \beta, A) : \gamma \in (-\varepsilon_0(K_0 A(\tau)), m_1),$$
$$\beta \in B_0, A \in \mathcal{A}^\circ_{(\gamma)}, 1/M < a < M\},$$

where $1 < M < \infty$, $\varepsilon_0(\cdot)$ is as defined in condition (C) and $K_0$ is as defined in Section 2.1. Let $D(\nu)$ be the space of all conditional densities $k(v|z)$ such that $k^e(v|z) L(v|z)^{1-e}$, where $L(v|z) = \int_v^\infty k(u|z)\,du$, is $\nu$-measurable. Denote $f_\psi(v|z) = p_\psi(v, 1|z)$. Also, for $f \in D(\nu)$, let $p_{(f)}(v, e|z) = f^e(v|z) S^{1-e}(v|z)$, where $S(v|z)$ is the survival function corresponding to $f(v|z)$. The main result of this section is as follows.

THEOREM 6. *Assume that conditions* (A1)–(A3), (B), (C), (D1)–(D3) *and* (H) *are satisfied by the data and the posited frailty distribution. Then for every $Q < \infty$ there exists an $\varepsilon > 0$ such that for each conditional density $f \in D(\nu)$, with $f \leq Q$ $\nu$-almost surely and $\int_\mathcal{X} |f - f_\psi|\,d\nu \leq \varepsilon$ for some $\psi \in \Psi_M$, there exists a unique Kullback–Leibler maximizer $\psi_{(f)} = \arg\max_{\psi_* \in \Psi_M} P_{(f)} \times \log(p_{\psi_*}/p_{(f)}) \in \Psi_M$, with $\sigma_{\psi_{(f)}} \equiv P_{(f)} \hat\sigma_{\psi_{(f)}} : H_\infty \mapsto H_\infty$ being one-to-one, where $P_{(f)} g = \int_\mathcal{X} g p_{(f)}\,d\nu$.*

REMARK 22. This theorem tells us that for any given class of proportional hazards frailty regression models parameterized by $\Psi_M$ and satisfying the stated regularity conditions, there exist an infinite number of true models not agreeing with the posited model but which satisfy conditions (F) and (G) with $\psi_* = (\gamma_*, \beta_*, A_*) = \psi_{(f)}$. In other words, conditions (F) and (G) are satisfied when the posited frailty family is sufficiently close to the true distribution. This is important in misspecified frailty model settings where uniqueness is not guaranteed by convexity. Note that without the condition bounding the true densities by $Q$, the Kullback–Leibler discrepancy between the true and posited models may be unbounded even when the respective densities are quite close in $L_1(\nu)$.

REMARK 23. It is worth emphasizing that $\gamma_*$ may be less than 0. In practice, one might mistakenly assume that the model is correctly specified and constrain the maximization to be over the subset of $\Psi_M$ for which $\gamma \geq 0$. Denote the resulting maximizer $\tilde\psi_* = (\tilde\gamma_*, \tilde\beta_*, \tilde A_*)$ and assume it is unique. The results in Section 2 can then be redone for the estimator $\tilde\psi_n$ defined in Remark 21. In the general setting, $\tilde\psi_n$ will be uniformly consistent for $\tilde\psi_*$, and $\sqrt{n}(\tilde\psi_n - \tilde\psi_*)$ will have three possible limiting distributions: when $\tilde\gamma_* > 0$, the limiting distribution is a Gaussian process as given in Theorem 4; when $\tilde\gamma_* = 0$ and the $\gamma$ term in $U_0^\tau(\tilde\psi_*) = 0$, the limiting distribution is a mixture of two Gaussian processes similar to the mixture described in Remark 21; and when $\tilde\gamma_* = 0$ but the $\gamma$ term in $U_0^\tau(\tilde\psi_*) < 0$, the limiting distribution for $\sqrt{n}(\tilde{\gamma}_n - \tilde{\gamma}_*)$ is a point mass at 0 while the remaining components have the limiting distribution resulting from assuming the Cox model ($\gamma_* = 0$). It is not possible, under the stated regularity conditions, to have $\tilde{\gamma}_* = 0$ but the $\gamma$ term in $U_0^\tau(\tilde{\psi}_*) > 0$, since this would imply that $\gamma_* > 0$. It also can be shown that the conditional limit law of the bootstrap, which imitates the foregoing estimation procedure, is equal to the limit law of $\sqrt{n}(\tilde{\psi}_n - \psi_*)$ in the sense of Corollary 1. Hence, when the bootstrap distribution of $\gamma_*$ under this constraint is frozen at 0, there is significant evidence against the frailty model being correctly specified.

6.2. *Identifying an independent covariate effect under a misspecified model.* In this section, we examine the effect of testing for a univariate covariate effect in the presence of other covariates and frailties when the posited frailty distribution and some of the covariates may be misspecified. We assume throughout this section that the data satisfy conditions (A1′), (A2) and (A3′) for the covariate process $Z = (Z_1, Z_2)$ and that the posited model satisfies conditions (B), (C), (D1), (D2) and (E1). We allow $Z_2$ to be misspecified and $Z_1$ to be partly misspecified, in that we only assume $P_0\{T \leq (\cdot) | Z_1 = z_1\}$ is monotone in $z_1$.

The following is the main result of this section.

THEOREM 7. *Assume conditions* (A1′), (A2), (A3′), (B), (C), (D1)–(D3) *and* (E1) *hold for the data and posited frailty model. Denote* $F_0^{z_1}(t) = P[T > t | Z_1 = z_1]$ *and assume that* $F_0^{z_1}(\cdot)$ *is monotone in* $z_1$ *almost surely. The covariates may be otherwise misspecified. Also assume condition* (F) *holds with Kullback–Leibler maximizer* $\psi_* = (\gamma_*, \beta_* = (\beta_{*1}, \beta_{*2}), A_*) \in \Psi_0$, *where* $\gamma_* \geq 0$. *Then* (i) *if* $F_0^{z_1}$ *is constant in* $z_1$, $\beta_{*1} = 0$; (ii) *if* $F_0^{z_1}$ *is strictly increasing in* $z_1$, $\beta_{*1} > 0$; *and* (iii) *if* $F_0^{z_1}$ *is strictly decreasing in* $z_1$, $\beta_{*1} < 0$.

REMARK 24. If we interpret the covariate effect of $Z_1$ to be positive when $F_0^{z_1}$ is increasing in $z_1$, negative when decreasing, 0 when constant and ambiguous otherwise, then Theorem 7 implies that the covariate effect can be consistently estimated up to the correct sign even if both $Z_1$ and $Z_2$ are otherwise misspecified. If condition (G) also holds, then the score and Wald tests for $H_0: \beta_{*1} = 0$ will be valid for testing the covariate effect of $Z_1$. These results generalize Kong and Slud (1997) to more general misspecification when fitting the more general class of models (1.1).

REMARK 25. Note that we require $\gamma_* \geq 0$. This is because condition (E1) appears to be needed for Theorem 7, and this condition only works when $\gamma \geq 0$. This requirement is stronger than necessary for the consistency and asymptotic normality results for possibly misspecified models given in Section 4.

6.3. *Coefficient effects under misspecified models.* In this section we examine the effect of regression parameter estimates under a misspecified frailty distribution and stronger conditions on the covariates. We assume that the data and posited frailty distribution satisfy (A1′), (A2), (A3″), (B), (C), (D1)–(D3) and (E1).

REMARK 26. Conditional linearity [condition (A3″)] is also used by Li and Duan (1989) in their study of regression analysis under link violation. Their results apply to fitting parametric models based on maximum likelihood and to semiparametric Cox models based on partial likelihood without censoring. Brillinger (1983) used this assumption to study unobserved Gaussian regressor variables in generalized linear models. In contrast, our results apply to semiparametric frailty regression models under frailty misspecification using nonparametric maximum likelihood with or without censoring. While condition (A3″) is sufficient, it may not be necessary for the results below.

We have the following proposition which extends Li and Duan (1989) to censoring when the true model has the form (1.1).

PROPOSITION 5. *Assume conditions* (A1′), (A2) *and* (A3″) *hold for the data, the posited model is a Cox proportional hazards model with* $\beta \in \mathbb{R}^d$, *and the true failure time distribution satisfies a proportional hazards frailty regression model with parameter* $\psi_0 = (\gamma_0, \beta_0, A_0) \in \Psi_0$, *where* $\gamma_0 \geq 0$ *and where the corresponding true negative log frailty transform family (denoted* $\{G_\gamma^\circ\}$*) satisfies conditions* (C), (D1)–(D3) *and* (E1). *Then conditions* (F) *and* (G) *are satisfied with* $\beta_* = \alpha_1 \beta_0$, *where* $\alpha_1 = 1$ *if* $\gamma_0 = 0$ *and* $\alpha_1 \in (0,1)$ *if* $\gamma_0 > 0$.

The following result establishes consistency up to scale when fitting (1.1) with a misspecified Laplace transform.

PROPOSITION 6. *Assume conditions* (A1′), (A2) *and* (A3″) *hold for the data; the posited and true proportional hazards frailty regression models satisfy conditions* (B), (C), (D1)–(D3) *and* (E1) *for the common index set* $\Psi$, *where the posited negative log frailty transform family is denoted* $\{G_\gamma\}$ *and the true family is denoted* $\{G_\gamma^\circ\}$; *and the true parameter value for the true model is* $\psi_0 = (\gamma_0, \beta_0, A_0) \in \Psi_0$, *where* $\gamma_0 \geq 0$. *Also assume condition* (F) *holds with Kullback–Leibler maximizer* $\psi_* = (\gamma_*, \beta_*, A_*) \in \Psi_0$, *where* $\gamma_* \geq 0$. *Then* $\beta_* = \alpha_2 \beta_0$, *where* $\alpha_2 = \alpha_1$ *when* $\gamma_* = 0$, $\alpha_2 > 0$ *when* $\gamma_* > 0$ *and* $\alpha_1$ *is as defined in Proposition* 5.

REMARK 27. If conditional linearity is violated, it may be possible to perform a reweighted maximum likelihood estimation procedure, based on weights described in Cook and Nachtsheim (1994). The resulting estimator would have the properties described in Proposition 6 and could be employed as a diagnostic for (A3″) by comparing to the unweighted estimator.

REMARK 28. Under the stated regularity conditions, these results show that using the Cox model when an unobserved frailty is present results in an estimate which is an attenuation of the true effect. When fitting (1.1) with the Laplace transform correctly specified, there is a deattenuation relative to the Cox model as a consequence of model identifiability. With a misspecified frailty distribution, the correct direction is obtained. However, it is unclear whether the effect size is deattenuated relative to the Cox model.

REMARK 29. One can test whether the Cox model is an attenuation of the true effect, that is, $\alpha_1 < 1$, if the score test for $H_0: \gamma_0 = 0$ remains valid and consistent under misspecification of the frailty distribution. Proving this in generality appears to be quite difficult, but the following result is a step in the right direction.

PROPOSITION 7. *Assume conditions* (A1′), (A2) *and* (A3″) *hold for the data and the posited and true proportional hazards frailty regression models satisfy conditions* (B), (C), (D1)–(D3), (E1) *and* (E2) *for the common index set* $\Psi$. *Then the score test for* $H_0: \gamma_0 = 0$ *based on the posited model is valid and consistent under positive contiguous alternatives.*

REMARK 30. Proposition 7 points out that this score test has the same form at $\gamma = 0$ under both correctly and incorrectly specified frailty models. Thus, to ensure that this score test is consistent for the fixed alternative $H_1: \gamma_0 > 0$, one would need to establish that the profile likelihood for $\gamma$, profiling over $\beta$ and $A$, is convex over the region $[0, m_1]$ when the model is correctly specified. This appears to be very challenging analytically.

**7. Computational issues.** We implemented the profile likelihood estimation method of Section 3.1, along with the bootstrap procedure of Section 4.4 based on Dirichlet weights, for the gamma and inverse Gaussian frailty models. We did not implement the log-normal frailty model because of the additional computational burden resulting from $\Lambda_\gamma$ not having a closed form. Although this issue can be addressed using Monte Carlo quadrature methods, we do not pursue it further here. To estimate the parameters in the gamma and inverse Gaussian frailty models, we maximized the nonparametric likelihood via profiling. A simple random search method based on the Metropolis–Hastings algorithm was used to maximize $pL_n(\theta)$ over $\Theta$.

For each candidate value of $\theta$, the fixed-point equation given in (3.2) was iterated until stabilization to obtain $\hat{A}_\theta$. Some simplification of $J_n^\psi$ occurs for these two frailty models. For the gamma frailty,

$$J_n^\psi(t) = \mathbb{P}_n\left[\frac{Y(t)e^{\beta'Z(t)}(1+\gamma\delta)}{1+\gamma H^\psi(V)}\right],$$

and for the inverse Gaussian frailty,

$$J_n^\psi(t) = \mathbb{P}_n\left[\frac{Y(t)e^{\beta'Z(t)}(\{1+2\gamma H^\psi(V)\}^{1/2}+\gamma\delta)}{1+\gamma H^\psi(V)}\right].$$

When the candidate value of $\gamma$ was negative, the likelihood was considered 0 if $G_\gamma(H^{(\theta,\hat{A}_\theta)}(V))$ was either negative or undefined for any data point. Overall, we found that this procedure was accurate and computationally efficient at finding the maximum, with or without bootstrap weights.

**8. Example: non-Hodgkin's lymphoma data.** The data are a subset of 1385 patients with aggressive non-Hodgkin's lymphoma (NHL), from 16 institutions and cooperative groups in North America and Europe. These patients were treated with a particular chemotherapy regimen. Survival was documented from start of treatment until either death or loss to follow-up. The censoring rate was 54.7%. Information on the following pretreatment covariates is complete for all patients in the subset: *age* at the diagnosis of NHL ($\leq 60$ or $> 60$ years), performance *status* (ambulatory or nonambulatory), serum lactate dehydrogenase *level* (below normal or above normal), number of extranodal disease *sites* ($\leq 1$ or $> 1$) and Ann Arbor classification of tumor *stage* [stage I or II (localized disease) or stage III or IV (advanced disease)]. Each characteristic is coded 0 for the first group in the parentheses and 1 for the second. These dichotomous predictors are the basis for the original model [Non-Hodgkin's Lymphoma Prognostic Factors Project (1993)]. A clinical reason for using dichotomous predictors is that it provides a simple classification of risk based on only a finite set of risk groups.

We now illustrate the utility of the procedures described in Section 7 for the gamma and inverse Gaussian frailty models. The bootstrap procedure based on the Dirichlet weights was repeated 500 times for inference. Parameter estimates, standard errors and $Z$ values for the parameters in the gamma frailty model both with unknown $\gamma$ (GF) and with the estimated value of $\gamma$ treated as known (GF$_0$), the inverse Gaussian frailty model with unknown $\gamma$ (IGF), the proportional odds model and the Cox model are given in Table 1. While the coefficient estimates are the same for GF and GF$_0$, the difference is that the standard errors for GF$_0$ are based on bootstrap estimates from a model with fixed $\gamma = 2.197$. The GF$_0$ standard errors are helpful in assessing the bias in precision estimation due to assuming $\gamma$ known. This bias is generally nontrivial and should not be ignored in practice.

TABLE 1
*Parameter estimates for the non-Hodgkin's lymphoma data using the gamma frailty model with $\gamma$ both unknown (GF) and fixed at the estimated value ($GF_0$), the inverse Gaussian frailty model with $\gamma$ unknown (IGF), the proportional odds model (PO) and the Cox model (PH)*

| Covariate | Parameter | Model | Estimate | S.E. | $Z$ value |
|---|---|---|---|---|---|
| (oddsrate) | $\gamma$ | GF | 2.197 | 0.447 | 4.914 |
| | | $GF_0$ | 2.197* | — | — |
| | | IGF | 4.325 | 1.855 | 2.332 |
| | | PO | 1* | — | — |
| | | PH | 0* | — | — |
| Age | $\beta_1$ | GF | 1.100 | 0.139 | 7.903 |
| | | $GF_0$ | 1.100 | 0.127 | 8.628 |
| | | IGF | 1.008 | 0.138 | 7.291 |
| | | PO | 0.881 | 0.118 | 7.457 |
| | | PH | 0.683 | 0.088 | 7.796 |
| Level | $\beta_2$ | GF | 1.024 | 0.159 | 6.433 |
| | | $GF_0$ | 1.024 | 0.144 | 7.155 |
| | | IGF | 0.933 | 0.143 | 6.524 |
| | | PO | 0.833 | 0.120 | 6.967 |
| | | PH | 0.624 | 0.092 | 6.796 |
| Status | $\beta_3$ | GF | 1.291 | 0.210 | 6.136 |
| | | $GF_0$ | 1.291 | 0.166 | 7.779 |
| | | IGF | 0.994 | 0.147 | 6.767 |
| | | PO | 0.949 | 0.145 | 6.562 |
| | | PH | 0.586 | 0.098 | 5.958 |
| Sites | $\beta_4$ | GF | 0.694 | 0.156 | 4.444 |
| | | $GF_0$ | 0.694 | 0.149 | 4.644 |
| | | IGF | 0.622 | 0.139 | 4.464 |
| | | PO | 0.546 | 0.122 | 4.458 |
| | | PH | 0.394 | 0.092 | 4.272 |
| Stage | $\beta_5$ | GF | 0.584 | 0.158 | 3.693 |
| | | $GF_0$ | 0.584 | 0.154 | 3.779 |
| | | IGF | 0.545 | 0.148 | 3.680 |
| | | PO | 0.485 | 0.137 | 3.549 |
| | | PH | 0.369 | 0.104 | 3.560 |

*$\gamma$ is fixed at the given value.

The attenuation of the covariate effects in the Cox model, predicted in Section 6.3, is evident in the results, although the attenuation does not appear to be uniform across all covariates. For *status*, the ratio of the parameter coefficient under the gamma frailty model to the coefficient under the Cox model is about 2.2, while the corresponding ratio for *stage* is only about 1.6. This difference may be related to the fact that the ratio of the standard errors of the *status* coefficient estimates for GF to $GF_0$ is about 1.27, while the corresponding ratio for *stage* is only 1.03. An anonymous referee has suggested that the Cox attenuation phenomenon for a covariate effect may depend on the degree to which that covariate's parameter estimate is correlated with the frailty variance estimate. Except for the $GF_0$ results, the $Z$ values for the covariate effects are fairly stable across models.

The estimated frailty variances in GF and IGF are 2.197 and 4.325, respectively, which are significantly higher than that assumed by the Cox model ($\gamma = 0$) and the proportional odds model ($\gamma = 1$). The maximized log profile likelihood values for the GF, IGF, PO and PH models are $-4618.39$, $-4628.37$, $-4623.30$ and $-4688.40$, respectively. This suggests that GF provides the best fit to the data and that $\gamma$ is significantly greater than 1 ($p = 0.0074$ via the two-sided Wald test based on the bootstrap). Also, PO is better that IGF, even though IGF is more flexible, seemingly.

In Figure 1, we plot the Kaplan–Meier estimates of the marginal survival distributions for the LDH level and performance status groups. The estimates from GF and the Cox model are also displayed. The survival estimate in a group (e.g., patients with *status* $= 0$) from GF is $\prod_{0 < s \leq t}\{1 - \Delta \bar{H}(s)\}$, where

$$\bar{H}(t) = \int_0^t \frac{\sum_i Y_i(s) \exp\{\hat{\beta}'_n Z_i(s)\}}{\sum_i Y_i(s)}$$
$$\times \left(1 + \hat{\gamma}_n \int_0^s \exp\{\hat{\beta}'_n Z_i(u)\} d\hat{A}_n(u)\right)^{-1} d\hat{A}_n(s),$$

and summation is over all observations in the group. That is, $\bar{H}(t)$ is a model-based estimate of the cumulative hazard which averages over the observed covariate distribution in that group. The estimate based on the Cox model uses the partial likelihood estimator, Breslow's estimator and 0 in the place of $\hat{\beta}_n$, $\hat{A}_n$ and $\hat{\gamma}_n$, respectively, in $\bar{H}$. The reason the Kaplan–Meier curves may be quite different from the model-based curves in the tail is that there are fewer observations in the subgroups available for the Kaplan–Meier curves, whereas the model-based curves utilize all of the data. In general, the GF estimates are closer to the Kaplan–Meier curves than the proportional hazards fit, particularly with the performance status group comparison. This demonstrates the superior fit of the frailty model. We also examined the survival curve estimates based on the proportional odds model and found them to be intermediate between GF and the Cox model. These are omitted from the figures for clarity.

Next, we illustrate the robust inference procedure for the best fitting survival probabilities under the assumed gamma frailty model. Survival predictions for two covariate values, representing an elderly high-risk patient $[Z = (1, 1, 1, 1, 1)']$ and an elderly low-risk patient $[Z = (1, 0, 0, 0, 0)']$, are shown in Figure 2. Also shown are 95% simultaneous confidence bands for the GF prediction using 500 multiplier bootstrap samples with Dirichlet weights. The Cox proportional hazards survival predictions are included for comparison. The Cox prediction for the high-risk patient significantly underestimates the long-term survival probability relative to GF. The difference between the Cox and GF predictions is less pronounced for the low-risk patient. Some improvement in the model fit may be possible if continuous rather than dichotomous covariates are used, but we do not pursue this further here.

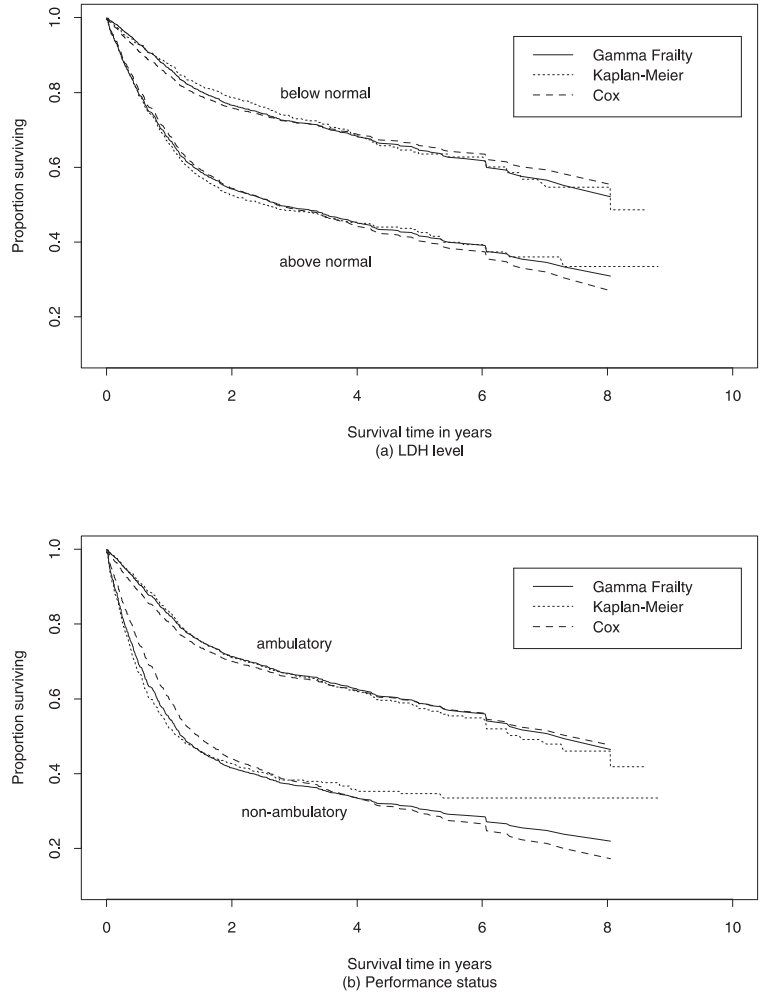

Fig. 1. *Estimated marginal survival distributions of non-Hodgkin's lymphoma patients for LDH level* (a) *and performance status* (b) *under the gamma frailty model. The Kaplan–Meier estimates and Cox model estimates are included for comparison.*

**9. Proofs.**

PROOF OF PROPOSITION 1. For any IGG($\alpha$) frailty family, the conditions hold with $\varepsilon_0(t) = (2/3)(1 \vee t)^{-1}$, $c_1(\gamma) = 1/\gamma$ and $c_2(\gamma) = (1-\alpha)/\gamma$, where $a \vee b$ is the maximum of $a$ and $b$. Establishing these results is straightforward. For the log-normal, we now show that the conditions hold with $\varepsilon_0(t) = (1 \vee t)^{-4}/64$, $c_1(\gamma) = 1$ and $c_2(\gamma) = 1/\gamma$. Complex analysis is involved since $\sqrt{\gamma}$ is imaginary for $\gamma < 0$. However, the imaginary components of $\Lambda_\gamma$, $G_\gamma$ and their derivatives are all 0. Moreover, $\Lambda_\gamma(t)$ and its first two deriva-

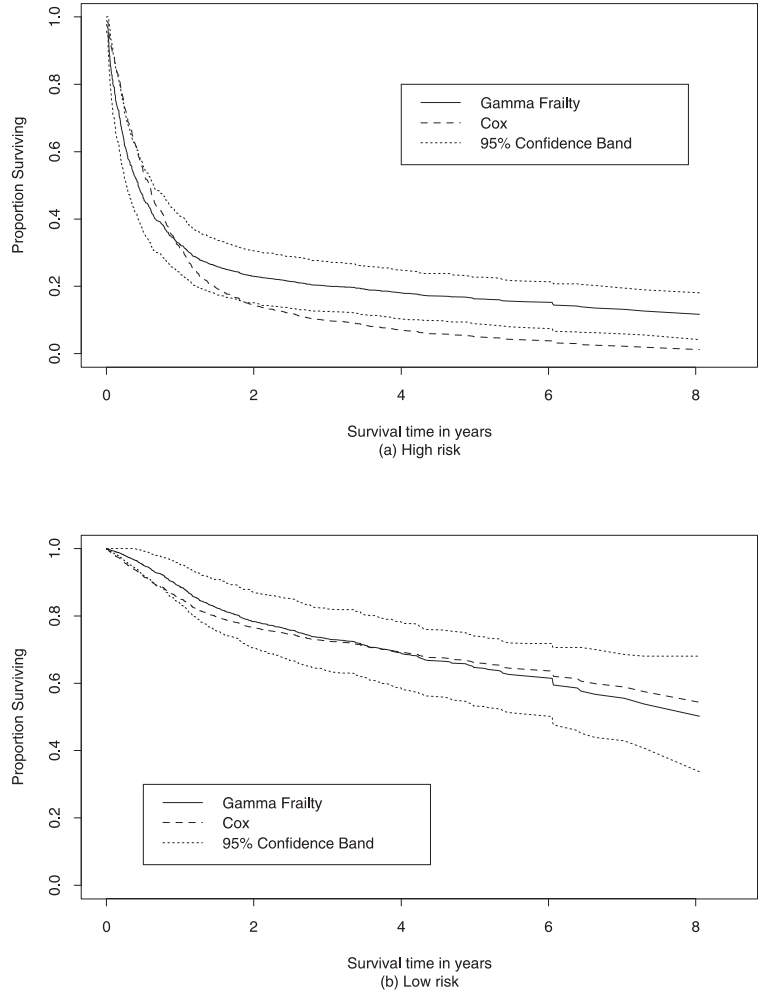

FIG. 2. *Survival predictions and 95% simultaneous confidence bands for the gamma frailty model: (a) high risk, $Z = (1,1,1,1,1)'$; (b) low risk, $Z = (1,0,0,0,0)'$. The Cox model predictions are included for comparison.*

tives in $t$ have the following form, with $\xi \equiv \sqrt{|\gamma|}$ and $u \equiv te^{\xi^2/2}$:

$$(9.1) \qquad (-1)^k e^{k\xi^2/2} \int_{\mathbb{R}} e^{-u\cos\xi v} \cos(u\sin\xi v - k\xi v)\phi(v)\,dv$$

for $k = 0, 1, 2$, respectively. If we establish that (9.1), for $k = 2$, is greater than 0 over the correct range, then (C) follows and showing (D1) is easy. Fix $0 < u' < \infty$. If there exists a $v_0 \geq 2$ and $\xi_0 \in [0, \pi/(2v_0)]$ such that

$$(9.2) \qquad u'\sin\xi_0 v_0 + 2\xi_0 v_0 = \pi/4$$

and such that the part of the integral over $|v| > v_0$ is completely dominated by the part over $|v| \leq v_0$, then (9.1) will be greater than 0 for all $u \in [0, u']$ and all $\xi \in [0, \xi_0]$. Assume that $v_0 \geq 2$ and $\xi_0$ satisfies (9.2). Then

$$\int_{-v_0}^{v_0} e^{-u'\cos\xi_0 v} \cos(u'\sin\xi_0 v - 2\xi_0 v)\phi(v)\,dv \geq \frac{0.95}{\sqrt{2}} e^{-u'}$$

and

$$\int_{|v|>v_0} e^{-u'\cos\xi_0 v} \cos(u'\sin\xi_0 v - 2\xi_0 v)\phi(v)\,dv \leq \frac{1}{\sqrt{2\pi}} e^{u' - v_0^2/2}.$$

Thus, the total integral is clearly positive when

$$\frac{1}{\sqrt{2\pi}} e^{u' - v_0^2/2} \left\{ \frac{0.95}{\sqrt{2}} e^{-u'} \right\}^{-1} \leq \frac{3}{4}.$$

Note that this is satisfied whenever $v_0 \geq 2\sqrt{1 \vee u'}$. Choosing any positive $\xi_0 \leq \pi(1 \vee u')^{-3/2}/24$ assures that there exists a $v_0 \geq 2\sqrt{1 \vee u'}$ which also satisfies (9.2). Setting $te^{\xi_0^2/2} = u'$, we have that $\xi_0 = (1/8)(1 \vee t)^{-3/2}$ is sufficient and then $\varepsilon_0(t) = (1 \vee t)^{-3}/64$ works. However, to satisfy (D2) and (D3) we reduce the rate to $(1 \vee t)^{-4}/64$. Lemma 3 gives that this rate is sufficient. □

LEMMA 3. *Conditions* (D2) *and* (D3) *are satisfied by the log-normal frailty model with* $\varepsilon_0(t) = (1 \vee t)^{-4}/64$, $c_1(\gamma) = 1$ *and* $c_2(\gamma) = 1/\gamma$.

PROOF. For (D2), let $W_k \equiv e^{\gamma_k^{1/2} Z - \gamma_k/2}$, where $Z \sim N(0,1)$. Now,

$$\sup_{u \geq 0} u\Lambda_{\gamma_k}(u) = \sup_{u \geq 0} \mathrm{E}[ue^{-uW_k}]$$

$$\leq \mathrm{E}\left[\sup_{u \geq 0} ue^{-uW_k}\right] \leq \mathrm{E}[W_k^{-1}],$$

$$\sup_{u \geq 0} |u^2 \dot{\Lambda}_{\gamma_k}(u)| = \sup_{u \geq 0} \mathrm{E}[u^2 W_k e^{-uW_k}]$$

$$\leq \mathrm{E}\left[\sup_{u \geq 0} u^2 W_k e^{-uW_k}\right] \leq \mathrm{E}[W_k^{-1}],$$

and it follows since $\mathrm{E}[W_k^{-1}] = e^{\gamma_k}$. For (D3), if $t_k \to \infty$ and $\gamma_k \to 0$, then $\gamma_k$ consists of one or both of two subsequences, one less than or equal to 0 and one greater than or equal to 0. Without loss of generality, assume $\gamma_k \to 0$ from above or below but not both. We begin with a sequence approaching from below. Let $\xi \equiv \sqrt{|\gamma|}$, $\xi_k \equiv \sqrt{|\gamma_k|}$, $u \equiv te^{\xi^2/2}$, $u_k \equiv t_k e^{\xi_k^2/2}$, and reparameterize $t\inf_{w \in [0,t]} \dot{G}_\gamma(w)$ for $\gamma < 0$ as

$$t \inf_{w \in [0,t]} \dot{G}_{-\xi^2}(w) = u \inf_{w \in [0,u]} \frac{\int_{\mathbb{R}} e^{-w\cos\xi v} \cos(w\sin\xi v - \xi v)\phi(v)\,dv}{\int_{\mathbb{R}} e^{-w\cos\xi v} \cos(w\sin\xi v)\phi(v)\,dv}$$

$$\equiv u \inf_{w \in [0,u]} g_\xi(w).$$

For $v > v_k \equiv u_k^{2/3}$, the $v^2$ term in $\phi(v)$ completely dominates $u_k$ since $v_k^2/u_k \to \infty$. Since $u_k \xi_k v_k \le u_k^{-1/3} \to 0$,

$$\inf_{w \in [0,u_k]} g_{\xi_k}(w) = \inf_{w \in [0,u_k]} \frac{\int_{\mathbb{R}} e^{w(1-\cos\xi_k v)} \cos(w\sin\xi_k v - \xi_k v)\phi(v)\,dv}{\int_{\mathbb{R}} e^{w(1-\cos\xi_k v)} \cos(w\sin\xi_k v)\phi(v)\,dv}$$
$$\to 1.$$

Hence, $u_k \inf_{w \in [0,u_k]} g_{\xi_k}(w) \to \infty$.

Now assume that $\gamma_k \to 0$ from above. Then

$$\frac{-t\dot{\Lambda}_\gamma(t)}{\Lambda_\gamma(t)} = \frac{\int_{\mathbb{R}} te^{\sqrt{\gamma}v - \gamma/2} \exp\{-te^{\sqrt{\gamma}v - \gamma/2}\}\phi(v)\,dv}{\int_{\mathbb{R}} \exp\{-te^{\sqrt{\gamma}v - \gamma/2}\}\phi(v)\,dv}$$

$$= \left[u^{-1} \frac{\int_{\mathbb{R}} e^{-w} \zeta_u(w) \phi_\gamma(w)\,dw}{\int_{\mathbb{R}} \zeta_u(w) \phi_\gamma(w)\,dw}\right]^{-1},$$

where $u \equiv te^{-\gamma/2}$, $\zeta_u(w) \equiv ue^w \exp\{-ue^w\}$ and $\phi_\gamma(w) \equiv \gamma^{-1/2}\phi(\gamma^{-1/2}w)$. Since $e^{-w}$ is a decreasing function, for any $w_k \to \infty$,

(9.3)
$$\frac{\int_{\mathbb{R}} e^{-w} \zeta_u(w) \phi_\gamma(w)\,dw}{\int_{\mathbb{R}} \zeta_u(w) \phi_\gamma(w)\,dw} \le \frac{\int_{-\infty}^{-w_k} e^{-w} \zeta_u(w) \phi_\gamma(w)\,dw}{\int_{-\infty}^{-w_k} \zeta_u(w) \phi_\gamma(w)\,dw}$$

$$= \frac{\int_{w_k}^{\infty} e^{w} \zeta_u(-w) \phi_\gamma(w)\,dw}{\int_{w_k}^{\infty} \zeta_u(-w) \phi_\gamma(w)\,dw}.$$

Now denote $u_k \equiv t_k e^{-\gamma_k/2}$ and let $w_k = \log(1 + \gamma_k u_k + u_k^{1/2})$. Since

$$\gamma_k^{-1} w_k (w_k + u_k e^{-w_k})^{-1} = w_k \left(\gamma_k w_k + \frac{\gamma_k u_k}{1 + \gamma_k u_k + u_k^{1/2}}\right)^{-1} \to \infty$$

for all $w \ge w_k$, the $\phi_\gamma(w)$ term dominates the expectation in (9.3) with $u = u_k$ and $\gamma = \gamma_k$. Hence, with this substitution, $(9.3) = 1 + \gamma_k u_k + u_k^{1/2} + O(1)$ and $t_k \dot{G}_{\gamma_k}(t_k) \to \infty$. The same arguments work when $\gamma_k \to \gamma > 0$, except that $\liminf_{k\to\infty} t_k \dot{G}_{\gamma_k}(t_k) \geq 1/\gamma$. Condition (D3) follows since (C) implies for $\gamma \geq 0$ that $\inf_{w \in [0,t]} \dot{G}_\gamma(w) = \dot{G}_\gamma(t)$. □

PROOF OF PROPOSITION 2. Fix the sample size $n$. If the conclusion of this proposition does not hold, there exists a sequence $\{\theta_m \equiv (\gamma_m, \beta_m) \in \overline{\Theta}\}$ so that, for $\tilde{A}_m \equiv \hat{A}_{\theta_m}$, either $\limsup_{m\to\infty} \tilde{A}_m(\tau) = \infty$ or $\liminf_{m\to\infty} \tilde{A}_m(\tau) = 0$. Assume first that $\limsup_{m\to\infty} \tilde{A}_m(\tau) = \infty$ and let $\{m_k\}$ be a subsequence for which $\lim_{k\to\infty} \tilde{A}_{m_k}(\tau) = \infty$ and $\gamma_{m_k} \to \gamma$. Let $\tilde{\psi}_m \equiv (\gamma_m, \beta_m, \tilde{A}_m)$ and $\psi_m \equiv (\gamma_m, \beta_m, \mathbb{P}_n N)$. Using arguments from the proof of Theorem 1, we can conclude that $\gamma > 0$ and, for a partition $0 = u_0 < u_1 < \cdots < u_J \leq \tau$, that

$$\begin{aligned}
\tilde{L}_n(\tilde{\psi}_m) &- \tilde{L}_n(\psi_m) \\
&\leq c + \log(\tilde{A}_{m_k}(\tau))\mathbb{P}_n[\delta \mathbf{1}\{V \in [u_{J-1}, \infty]\} \\
&\qquad - (c_1(\gamma) + \delta)\mathbf{1}\{V \in [u_J, \infty]\}] \\
&+ \sum_{j=1}^{J-1} \log(\tilde{A}_{m_k}(u_j))\mathbb{P}_n[\delta \mathbf{1}\{V \in [u_{j-1}, u_j]\} \\
&\qquad - (c_1(\gamma) + \delta)\mathbf{1}\{V \in [u_j, u_{j+1}]\}],
\end{aligned} \tag{9.4}$$

where $c \in (0, \infty)$ is a constant not depending on the parameter values or on the partition, and where the summation is 0 when $J = 1$. Let $T_1, \ldots, T_J$ be the observed failure times and define a partition with $u_J = T_J$, and, if $J > 1$, let $u_j \in (T_j, T_{j+1})$ for $j = 1, \ldots, J-1$. Now the intervals $[u_j, u_{j+1}]$, for $j = 0, \ldots, J$, where $u_{J+1} = \infty$, all contain exactly one failure time. Thus, $(9.4) \leq c - c_1(\gamma)[\log(\tilde{A}_{m_k}(\tau)) + \sum_{j=1}^{J-1} \log(\tilde{A}_{m_k}(u_j))]$. Hence, using again arguments from the proof of Theorem 1, $(9.4) \to -\infty$. This is a contradiction. Thus, $\limsup_{m\to\infty} \hat{A}_m(\tau) < \infty$. The proof that $\inf_{\theta \in \overline{\Theta}} \hat{A}_\theta(\tau) > 0$ can also be obtained from arguments similar to those employed in the proof of Theorem 1. □

PROOF OF THEOREM 1. Let $(\mathcal{X}^\infty, \mathcal{B}^\infty, P_0^\infty)$ be the probability space for infinite sequences of observations, let $W \subset \mathcal{X}^\infty$ be the set of observation sequences for which $\mathbb{P}_n N$ converges uniformly to $\mu_0 \equiv P_0 N$ and note that $P_*(W) = 1$. Then if the conclusion of Theorem 1 does not hold, there exist a sequence $\{\theta_n \equiv (\gamma_n, \beta_n) \in \overline{\Theta}\}$ and an $\omega \in W$ so that, for $\hat{A}_n \equiv \hat{A}_{\theta_n}\{\omega\}$ (we will suppress dependence on $\omega$ hereafter), $\limsup_{n\to\infty} \hat{A}_n(\tau) = \infty$, $\liminf_{n\to\infty} \hat{A}_n(\tau) = 0$ or $\hat{A}_n$ is not asymptotically close to an absolutely continuous function with bounded derivative.

Assume first that $\limsup_{n\to\infty} \hat{A}_n(\tau) = \infty$. Now suppose $\gamma_n$ has an accumulation point at $\gamma < 0$ along a subsequence for which $\hat{A}_n(\tau) \to \infty$. But this is impossible by (C) and the constraints on $\mathcal{A}_{(\gamma)}$. Thus, $\gamma_n$ has no such accumulation points less than 0. Suppose, however, that $\gamma_n$ has 0 as one of these accumulation points. Accordingly, take a subsequence $\{n_k\}$ such that $\gamma_{n_k} \to 0$ and $\hat{A}_{n_k}(\tau) \to \infty$. Since, by (D1), $\dot{G}_\gamma(u) - \ddot{G}_\gamma(u)\{\dot{G}_\gamma(u)\}^{-1} = -\ddot{\Lambda}_\gamma(u)\{\dot{\Lambda}_\gamma(u)\}^{-1} \geq 0$, we have by (3.3) that

$$\hat{A}_{n_k}(\tau) \leq \int_0^\tau \left(\mathbb{P}_n\left[Y(t)e^{\beta'_{n_k}Z(t)}\left(\dot{G}_{\gamma_{n_k}}\{H^{\hat{\psi}_{n_k}}(V)\} - \delta\frac{\ddot{G}_{\gamma_{n_k}}\{H^{\hat{\psi}_{n_k}}(V)\}}{\dot{G}_{\gamma_{n_k}}\{H^{\hat{\psi}_{n_k}}(V)\}}\right)\right]\right)^{-1}$$
$$\times d\mathbb{P}_n\{N(t)\}$$
$$\leq O(1)\left\{\inf_{u \in [0, K_0\hat{A}_{n_k}(\tau)]} \dot{G}_{\gamma_{n_k}}(u)\right\}^{-1},$$

since $\inf_{t \in [0,\tau]} \mathbb{P}_n Y(t)e^{\beta'_{n_k}Z(t)}$ is bounded below for all $n$ large enough. Thus, by (D3), $1 \leq O(1) \{\hat{A}_{n_k}(\tau)\inf_{u \in [0, K_0\hat{A}_{n_k}(\tau)]} \dot{G}_{\gamma_{n_k}}(u)\}^{-1} \to 0$, which is a contradiction. Hence, for any subsequence with $\hat{A}_{n_k}(\tau) \to \infty$, the accumulation points of $\gamma_{n_k}$ are greater than 0. Now let $\{n_k\}$ be a subsequence with $\hat{A}_{n_k}(\tau) \to \infty$ and $\gamma_{n_k} \to \gamma > 0$.

Let $\hat{\psi}_n \equiv (\gamma_n, \beta_n, \hat{A}_n)$ and $\psi_n \equiv (\gamma_n, \beta_n, \mathbb{P}_n N)$. Then

$$0 \leq \tilde{L}_{n_k}(\hat{\psi}_{n_k}) - \tilde{L}_{n_k}(\psi_{n_k})$$
$$\leq \mathbb{P}_{n_k}\left\{\int_0^\tau \left[\log\left\{\frac{\dot{G}_{\gamma_{n_k}}(H^{\hat{\psi}_{n_k}}(s))}{\dot{G}_{\gamma_{n_k}}(H^{\psi_{n_k}}(s))}\right\} + \log(n_k\Delta\hat{A}_{n_k}(s))\right] dN(s)\right.$$
$$\left. + \log\left[\frac{\Lambda_{\gamma_{n_k}}(H^{\hat{\psi}_{n_k}}(V))}{\Lambda_{\gamma_{n_k}}(H^{\psi_{n_k}}(V))}\right]\right\}$$
$$\leq O(1) + \mathbb{P}_{n_k}\left\{\int_0^\tau \log(n_k\Delta\hat{A}_{n_k}(s)) dN(s)\right.$$
$$+ [-\delta(1 + c_1(\gamma))\log H^{\hat{\psi}_{n_k}}(V)] \wedge 0$$
$$\left. - (1-\delta)G_{\gamma_{n_k}}(H^{\hat{\psi}_{n_k}}(V))\right\}$$
$$\equiv C_{n_k},$$

since

$$\mathbb{P}_n\left[\int_0^\tau \log\{\dot{G}_{\gamma_{n_k}}(H^{\psi_{n_k}}(s))\} dN(s) + \log\Lambda_{\gamma_{n_k}}(H^{\psi_{n_k}}(V))\right] = O(1)$$

and since (D2) implies

$$\log\dot{G}_{\gamma_{n_k}}(H^{\hat{\psi}_{n_k}}(s)) \leq [-(1 + c_1(\gamma))\log H^{\hat{\psi}_{n_k}}(s)] \wedge 0 + O(1).$$

Now $C_{n_k}$ is bounded above by

$$O(1) + \mathbb{P}_{n_k}\left\{\int_0^\tau \log(n_k \Delta \hat{A}_{n_k}(s))\,dN(s)\right.$$
(9.5)
$$\left. + [-(\delta + c_1(\gamma))\log \hat{A}_{n_k}(V)] \wedge 0\right\},$$

since (D2) also implies $G_{\gamma_{n_k}}(H^{\hat{\psi}_{n_k}}(V)) \geq [c_1(\gamma)\log \hat{A}_{n_k}(V)] \vee 0 + O(1)$. For a sequence $0 = u_0 < u_1 < u_2 < \cdots < u_J = \tau$, let $N^j(s) \equiv N(s)\mathbf{1}\{V \in [u_{j-1}, u_j]\}$, $j = 1,\ldots,J$. By Jensen's inequality,

$$\int_0^\tau \log(n_k \Delta \hat{A}_{n_k}(s))\,d\mathbb{P}_{n_k}\{N^j(s)\}$$
$$\leq \mathbb{P}_{n_k} N^j(\tau) \log\left(\int_0^{u_j} n\Delta \hat{A}_{n_k}(s)\,d\mathbb{P}_{n_k} N^j(s)/\mathbb{P}_{n_k} N^j(\tau)\right)$$
$$\leq O(1) + \log(\hat{A}_{n_k}(u_j))\mathbb{P}_{n_k}(\delta \mathbf{1}\{V \in [u_{j-1}, u_j]\}).$$

Thus, (9.5) is dominated by

$$O(1) + \log(\hat{A}_{n_k}(\tau))\mathbb{P}_{n_k}[\delta \mathbf{1}\{V \in [u_{J-1}, \infty]\} - (c_1(\gamma) + \delta)\mathbf{1}\{V \in [\tau, \infty]\}]$$
(9.6)
$$+ \sum_{j=1}^{J-1} \log(\hat{A}_{n_k}(u_j))\mathbb{P}_{n_k}[\delta \mathbf{1}\{V \in [u_{j-1}, u_j]\}$$
$$- (c_1(\gamma) + \delta)\mathbf{1}\{V \in [u_j, u_{j+1}]\}].$$

Choose $\varepsilon: 0 < \varepsilon < P_0\{V = \tau\}$ and the sequence $\{u_j\}$ for finite $J$ such that $c_1(\gamma)\varepsilon/(c_1(\gamma)+1) < \mu_0(u_{j+1}) - \mu_0(u_j) < c_1(\gamma)\varepsilon$ for $j = 0, \ldots, J-2$ and $\mu_0(\tau) - \mu_0(u_{J-1}) = c_1(\gamma)\varepsilon/(c_1(\gamma)+1)$. Note that (3.3) implies $\log(\hat{A}_{n_k}(u)) \geq \log(\mathbb{P}_{n_k} N(u)) + O(1)$, since $\dot{G}_{\gamma_{n_k}}(t) - \delta \ddot{G}_{\gamma_{n_k}}(t)/\dot{G}_{\gamma_{n_k}}(t) \leq 2 - \ddot{G}_{\gamma_{n_k}}(0+)$ for all $t \geq 0$ and all $k$ large enough. Since also $\mathbb{P}_{n_k} N \to \mu_0$ uniformly, (9.6) goes to $-\infty$. This is a contradiction. Thus, $\limsup \hat{A}_{n_k}(\tau) < \infty$.

Now assume that there are a sequence $\{(\gamma_n, \beta_n) \in \overline{\Theta}\}$ and an $\omega \in W$ so that $\liminf_{n \to \infty} \hat{A}_n(\tau) = 0$. Define $\tilde{A}_n = \varepsilon_0^{-1}(c_0)\mathbb{P}_n N/K_0$ and note that $\tilde{A}_n \in \mathcal{A}_{(\gamma_n)}$ for all $n \geq 1$. Now let $\{n_k\}$ be a subsequence with $\hat{A}_{n_k}(\tau) \to 0$ and define $\psi_n \equiv (\gamma_n, \beta_n, \tilde{A}_n)$. Then

$$0 \leq \tilde{L}_{n_k}(\hat{\psi}_{n_k}) - \tilde{L}_{n_k}(\psi_{n_k})$$
$$\leq O(1) + \mathbb{P}_{n_k}\left\{\int_0^\tau \log(\Delta \hat{A}_{n_k}(s))\,dN(s)\right\}$$
$$\leq O(1) + \mathbb{P}_{n_k}\left\{\int_0^\tau \log(\hat{A}_{n_k}(\tau))\,dN(s)\right\} \to -\infty.$$

This is again a contradiction. Hence, $\liminf \hat{A}_n(\tau) > 0$. By previous arguments and (3.3), we also have, for $s, t \in [0, \tau]$, that

$$|\hat{A}_n(s) - \hat{A}_n(t)| \leq O(1)\mathbb{P}_n|N(s) - N(t)| \leq O(1)|\mu_0(s) - \mu_0(t)| + o(1),$$

and the conclusions of the theorem hold by (A2). $\square$

PROOF OF LEMMA 1. Fix a convergent sequence $\{\theta_m\} \in \overline{\Theta}$. The arguments used in the proof of Theorem 1, after replacing the measure $\mathbb{P}_n$ with $P_0$, can be used with only minor modification to show that $\limsup_{m\to\infty} A_{\theta_m}(\tau) < \infty$ and $\liminf_{m\to\infty} A_{\theta_m}(\tau) > 0$. Using again arguments from the proof of Theorem 1, we can also establish that $|A_{\theta_m}(s) - A_{\theta_m}(t)| \leq c|\mu_0(s) - \mu_0(t)|$ for all $s, t \in [0, \tau]$ and a constant $c \in (0, \infty)$ not depending on the sequence. The desired results now follow from condition (A2). $\square$

PROOF OF THEOREM 2. Let $W \subset \mathcal{X}^\infty$ be the set of data sequences for which $\mathbb{P}_n N \to \mu_0$ uniformly and note that the class of functions

$$\mathcal{G}_k \equiv \left\{ Y(t)e^{\beta' Z(t)} \left( \dot{G}_\gamma\{H^\psi(V)\} - \delta \frac{\ddot{G}_\gamma\{H^\psi(V)\}}{\dot{G}_\gamma\{H^\psi(V)\}} \right) : \right.$$
$$\left. t \in [0, \tau], \psi \in \Psi \text{ and } A(\tau) \leq k \right\}$$

is $P_0$-Glivenko–Cantelli for each $k < \infty$. To see this, arguments given in the proof of Proposition 8 verify that the classes $\{Y(t)e^{\beta' Z(t)} : t \in [0, \tau], \beta \in \overline{B}_0\}$ and $\{H^\psi(V) : \psi \in \Psi, A(\tau) \leq k\}$ are Donsker; conditions (C) and (D1) imply that the maps $(\gamma, t) \mapsto \dot{G}_\gamma(t)$, $(\gamma, t) \mapsto \ddot{G}_\gamma(t)$ and $(\gamma, t) \mapsto [\dot{G}_\gamma(t)]^{-1}$ are bounded and Lipschitz over the domain $[-\varepsilon_0(u), m] \times [0, u]$ for any $u \in (0, \infty)$; and $\Lambda(\tau) \leq k$ implies that $H^\psi(V) \leq kK_0$ for all $\psi \in \Psi$ almost surely. Thus, the classes $\{\dot{G}_\gamma\{H^\psi(V)\} : \psi \in \Psi, A(\tau) \leq k\}$, $\{\ddot{G}_\gamma\{H^\psi(V)\} : \psi \in \Psi, A(\tau) \leq k\}$ and $\{[\dot{G}_\gamma\{H^\psi(V)\}]^{-1} : \psi \in \Psi, A(\tau) \leq k\}$ are all Donsker by Theorem 2.10.6 of VW. Since products of bounded Donsker classes are Donsker, the class $\mathcal{G}_k$ is Donsker and hence also Glivenko–Cantelli.

For each $M < \infty$, let $W_M \subseteq W$ be the subset of data sequences for which the limit points of $\{\hat{A}_\theta : \theta \in \overline{\Theta}\}$ are in $\mathcal{K}_M$ and also for which $\mathbb{P}_n \to P_0$ in $\ell^\infty(\mathcal{G}_M)$. For any $\theta \in \overline{\Theta}$, let $\tilde{A}_\theta(t) \equiv \int_0^t [J_0^{\psi_\theta}(s) + \rho_0(\psi_\theta)]^{-1}\mathbb{P}_n\{dN(s)\}$, where $J_0^\psi \equiv P_0 J_n^\psi$ and

$$\rho_0(\psi) \equiv \frac{P_0 \delta - \int_0^\tau J_0^\psi(s)\, dA(s)}{A(\tau)}.$$

Also let $\hat{\psi}_\theta = (\theta, \hat{A}_\theta)$ and $\tilde{\psi}_\theta \equiv (\theta, \tilde{A}_\theta)$. We will first show that

(9.7) $$\sup_{\theta \in \overline{\Theta}} |\tilde{L}_n(\hat{\psi}_\theta) - \tilde{L}_n(\tilde{\psi}_\theta)| \to 0$$

outer almost surely and then show that

$$\sup_{\theta_1,\theta_2\in\overline{\Theta}}|\tilde{L}_n(\tilde{\psi}_{\theta_1}) - \tilde{L}_n(\tilde{\psi}_{\theta_2}) - I_0(\psi_{\theta_1},\psi_{\theta_2})| \to 0 \quad (9.8)$$

outer almost surely, and the proof will be complete.

Fix $M < \infty$, choose a $w \in W_M$ and let $\{n\}$ index the corresponding data sequence. Let $\{\theta_n\}$ be any parallel sequence of parameters in $\overline{\Theta}$ and let $\{n_k\}$ be any convergent subsequence with $\theta_{n_k} \to \theta^*$, $\hat{\psi}_{\theta_{n_k}} \to \psi^* = (\theta^*, A^*)$ and $\tilde{\psi}_{\theta_{n_k}} \to \psi^{**} = (\theta^*, A_{\theta^*})$. The last convergence statement follows from the definitions of $\tilde{A}_\theta$ and $A_\theta$. Since

$$\begin{aligned}(9.9)\qquad \frac{d\hat{A}_{n_k}(t)}{d\tilde{A}_{n_k}(t)} &= \frac{J_0^{\tilde{\psi}_{\theta_{n_k}}}(t) + \rho_0(\tilde{\psi}_{\theta_{n_k}})}{J_n^{\hat{\psi}_{\theta_{n_k}}}(t) + \rho_n(\hat{\psi}_{\theta_{n_k}})} \\ &\to \frac{J_0^{\psi^{**}}(t) + \rho_0(\psi^{**})}{J_0^{\psi^*}(t) + \rho_0(\psi^*)},\end{aligned}$$

uniformly over $t \in [0, \tau]$, where the limit of $(9.9) = dA^*(t)/dA_{\theta^*}(t)$ is bounded below and in total variation, we have that

$$\mathbb{P}_n \int_0^\tau \log\{d\hat{A}_{\theta_{n_k}}(t)/d\tilde{A}_{\theta_{n_k}}(t)\}\, dN(t) \to \int_0^\tau \log\{dA^*(t)/dA_{\theta^*}(t)\}\, d\mu_0(t).$$

Hence, it follows that

$$0 \leq \tilde{L}_{n_k}(\hat{\psi}_{\theta_{n_k}}) - \tilde{L}_{n_k}(\tilde{\psi}_{\theta_{n_k}}) \to P_0 \log \frac{p_{\psi^*}}{p_{\psi^{**}}} \leq 0.$$

Since this is true for every such convergent subsequence, and since $M < \infty$ can be increased so that $P_*(W_M)$ is arbitrarily close to 1, we have established (9.7).

Since by Lemma 1

$$\frac{dA_{\theta_1}(t)}{dA_{\theta_2}(t)} = \frac{J_0^{\psi_{\theta_2}}(t) + \rho_0(\psi_{\theta_2})}{J_0^{\psi_{\theta_1}}(t) + \rho_0(\psi_{\theta_1})}$$

is bounded below and in total variation, uniformly over $\theta_1, \theta_2 \in \overline{\Theta}$, we have that

$$\mathbb{P}_n \int_0^\tau \log\{dA_{\theta_1}(t)/dA_{\theta_2}(t)\}\, dN(t) \to \int_0^\tau \log\{dA_{\theta_1}(t)/dA_{\theta_2}(t)\}\, d\mu_0(t),$$

uniformly over $\theta_1, \theta_2 \in \overline{\Theta}$. Hence, it follows that (9.8) holds. □

PROOF OF LEMMA 2. By the smoothness assumed in (D1) of the involved derivatives,

$$\lim_{t\downarrow 0} \sup_{\psi\,:\,\|\psi\|_{(r)}\leq 1} \sup_{h\in H_r} \left|\int_0^1 \psi(\sigma_{\psi_1+st\psi}(h) - \sigma_{\psi_1}(h))\, ds\right| = 0.$$

Thus, $\sup_{h \in H_r} |U_0^\tau(\psi_1 + \psi)(h) - U_0^\tau(\psi_1)(h) + \psi(\sigma_{\psi_1}(h))| = o(\|\psi\|_{(r)})$ as $\|\psi\|_{(r)} \to 0$. $\square$

PROOF OF THEOREM 3. Recall that $\hat{A}_n \equiv \hat{A}_{\hat{\theta}_n}$ and $\hat{\psi}_n = (\hat{\theta}_n, \hat{A}_n)$, where $\hat{\theta}_n$ is the profile MLE. Theorem 1 implies that the set $W \subset \mathcal{X}^\infty$ of data sequences for which the limit points of $\hat{A}_n$ are in $\mathcal{K}_M$, for some $M < \infty$, has inner probability 1. Accordingly, fix the data sequence $w \in W$ and take a subsequence $\{n_k\}$ for which $\hat{\psi}_{n_k}$ converges uniformly to some $\psi = (\theta, A) \in \Psi$, with $A \in \mathcal{K}_M$ for some $M < \infty$. Let $\psi_n = (\theta_*, \tilde{A}_{\theta_*})$, where $\theta_* \equiv (\gamma_*, \beta_*)$ and $\tilde{A}_\theta$ is as defined in the proof of Theorem 2. By Theorem 2 we have

$$0 \leq \tilde{L}_{n_k}(\hat{\psi}_{n_k}) - \tilde{L}_{n_k}(\psi_{n_k}) \to P_0 \log \frac{p_\psi}{p_{\psi_*}} \leq 0,$$

and hence condition (F) implies that $\psi = \psi_*$. Since this is true for every convergent subsequence, $\hat{\psi}_n \to \psi_*$ with inner probability 1. Since $\hat{A}_n$ is a piecewise constant function with mass $\Delta \hat{A}_n$ only at observed failure times $t_1, \ldots, t_{m_n}$, $\hat{\psi}_n$ is a functional of a maximum taken over $m_n + d + 1$ real variables. This structure implies that $\sup_{t \in [0,\tau]} |\hat{A}_n(t) - A_*(t)|$ is a measurable random variable, and hence the uniform distance between $\hat{\psi}_n$ and $\psi_*$ is also measurable. Thus, the convergence with inner probability 1 can be strengthened to outer almost sure convergence. $\square$

PROOF OF PROPOSITION 3. Since $f_0$ does not depend on covariates, there exists an $A_0 \in \mathcal{A}_0$ so that $S_0(t) = \exp(-A_0(t))$ for all $t \in [0, \tau]$. Hence, the parameter value $\psi_0 = (\gamma_0 = 0, \beta_0 = 0, A_0)$ for the posited model describes the true distribution of the failure times. Arguments in Theorem 3 now yield that, with inner probability 1, all limit points of the maximum likelihood estimator $\hat{\psi}_n$ lie in a compact set $\tilde{\Psi}$ for which $P_0 \log(p_\psi/p_0) = 0$. Since $P_0 \log(p_\psi/p_0) \leq -\int_{\mathcal{X}} (p_\psi^{1/2} - p_0^{1/2})^2 p_0 \, d\nu$ and $p_0 = p_{\psi_0}$, we now know that all Kullback–Leibler maximizers $\psi = (\gamma, \beta, A) \in \tilde{\Psi}$ must satisfy $G_\gamma(H^\psi(t)) = G_0(H^{\psi_0}(t))$ for all $t \in [0, \tau]$. This implies that $\beta = \beta_0$ by arguments given in Proposition 4. The desired result now follows. $\square$

PROOF OF THEOREM 4. From Section 3.3, we have $\sigma_{\psi_*}^{13}(h_3) = \int_0^\tau f_1(s) \times h_3(s) \, dA_*(s)$, $\sigma_{\psi_*}^{23}(h_3) = \int_0^\tau f_2(s) h_3(s) \, dA_*(s)$ and $\sigma_{\psi_*}^{33}(h_3) = g_1(\cdot) \int_0^\tau f_3(s) \times h_3(s) \, dA_*(s) + g_2(\cdot) h_3(\cdot)$, where $f_1, f_3, g_1 : \mathbb{R} \mapsto \mathbb{R}$ and $f_2 : \mathbb{R} \mapsto \mathbb{R}^d$ are bounded and where $g_2(s) \equiv P_0\{\hat{\xi}_{\psi_*}^{(0)} Y(s) e^{\beta_*' Z(s)}\}$. From the proof of Theorem 1, $0 < g_2(s) < \infty$ for all $s \in [0, \tau]$. Thus, $\sigma_{\psi_*} = \sigma_{\psi_*}^{(1)} + \sigma_{\psi_*}^{(2)}$, where

$$\sigma_{\psi_*}^{(1)}(h) \equiv \begin{pmatrix} 1 & 0 & 0 \\ 0 & I & 0 \\ 0 & 0 & g_2(\cdot) \end{pmatrix} \begin{pmatrix} h_1 \\ h_2 \\ h_3 \end{pmatrix}$$

is continuously invertible and onto, $\sigma_{\psi_*}^{(2)} = \sigma_{\psi_*} - \sigma_{\psi_*}^{(1)}$ is compact and $I$ denotes the identity. Since $\sigma_{\psi_*}$ is one-to-one by condition (G), $\sigma_{\psi_*}$ is continuously invertible and onto, with inverse $\sigma_{\psi_*}^{-1}$. This now implies that, for each $r > 0$, there is an $s > 0$ with $\sigma^{-1}(H_s) \subset H_r$. Fix $r > 0$. Continuous invertibility of $-\dot{U}_{\psi_*}$ on $\lin\Psi$, where lin denotes linear span, now follows by Proposition A.1.7 of BKRW since

$$\inf_{\psi \in \lin\Psi} \frac{\|\dot{U}_{\psi_*}(\psi)\|_{(r)}}{\|\psi\|_{(r)}} \geq \inf_{\psi \in \lin\Psi} \frac{\sup_{h \in \sigma_{\psi_*}^{-1}(H_q)} |\psi(\sigma_{\psi_*}(h))|}{\|\psi\|_{(r)}}$$

$$= \inf_{\psi \in \Psi} \frac{\|\psi\|_{(q)}}{\|\psi\|_{(r)}} \geq \frac{q}{3r}.$$

This proposition also implies that $-\dot{U}_{\psi_*}$ is onto. By Proposition 8,

$$(9.10) \quad \begin{aligned} &\sqrt{n}(U_n^\tau(\hat{\psi}_n)(h) - U_0^\tau(\hat{\psi}_n)(h)) - \sqrt{n}(U_n^\tau(\psi_*)(h) - U_0^\tau(\psi_*)(h)) \\ &= o_P(1 + \sqrt{n}\|\hat{\psi}_n - \psi_*\|), \end{aligned}$$

uniformly over $h \in H_r$, where $\|\cdot\|$ is the uniform metric.

Applying the $Z$-estimator master theorem (Theorem 3.3.1 of VW) now gives the desired weak convergence of $\sqrt{n}(\hat{\psi}_n - \psi_*)$, provided $U_0^\tau(\psi_*)(\cdot) = 0$, $U_n^\tau(\hat{\psi}_n)(\cdot) = 0$ asymptotically and $U^\tau(\psi_*)(\cdot)$ is $P_0$-Donsker. The first two conditions follow from condition (F) and the fact that $\hat{\psi}_n$ is asymptotically an interior point in $\Psi$ by consistency and condition (C). Because products of bounded Donsker classes are Donsker, showing $\delta h_3(V)$ and $\int_0^\tau Y(s) e^{\beta_*' Z(s)} h_3(s) \, dA_*(s)$ are Donsker, as processes indexed by $h_3: (h_1, h_2, h_3) \in H_r$, is sufficient, since $\beta_*$ and $A_*$ are fixed. First, since all functions in $H_r$ are bounded in total variation, $\delta h_3(V)$ is Donsker, as a class indexed by $h_3$, since it is the product of bounded Donsker classes. Next, $\{\beta_*' Z(t) : t \in [0, \tau]\}$ is Donsker since the total variation of $Z$ is less than or equal to $m_0$ with probability 1 by condition (A3). Since $\exp(\cdot)$ is Lipschitz on compacts and $\{Y(t), t \in [0, \tau]\}$ is monotone and bounded, $\{Y(t) e^{\beta_*' Z(t)}\}$ is Donsker. Finally, since the map from $Y(\cdot) e^{\beta_*' Z(\cdot)}$ to $\int_0^\tau Y(s) e^{\beta_*' Z(s)} h_3(s) \, dA_*(s)$, as a map from an element in $\ell^\infty([0, \tau])$ to $\ell^\infty(H_r)$, is continuous and linear, the continuous mapping theorem yields the desired Donsker property. These results now imply that $\hat{\psi}_n(h)$ is asymptotically linear with influence function $\tilde{\ell}(h) \equiv U^\tau(\psi_*)(\sigma_{\psi_*}^{-1}(h))$ and covariance $V_*(g, h) = \mathrm{E}[\tilde{\ell}(g)\tilde{\ell}(h)]$ for $g, h \in H_r$. Taking $r \geq 1$ yields weak convergence in the uniform metric, since $H_1$ is sufficiently rich as noted earlier. □

PROPOSITION 8. *Expression* (9.10) *holds.*

PROOF. If for some $\varepsilon > 0$ $\{U^\tau(\psi)(h) - U^\tau(\psi_*)(h) : \|\psi - \psi_*\| < \varepsilon, h \in H_r\}$ is $P_0$-Donsker and $\lim_{\psi \to \psi_*} \sup_{h \in H_r} P_0\{U^\tau(\psi)(h) - U^\tau(\psi_*)(h)\}^2 = 0$, then (9.10) holds by Lemma 3.3.5 of VW. The latter condition follows from condition (D1). The Donsker condition requires more work. Let $\Psi_\varepsilon \equiv \{\psi : \|\psi - \psi_*\| \leq \varepsilon\}$. Take $\varepsilon$ small enough so that $\Psi_\varepsilon \subset \Psi$: such an $\varepsilon$ always exists by (C) and the fact that $\gamma_* \geq 0$. Because $Z$ has bounded total variation and the class $\{\beta, \beta \in \overline{B}_0\}$ is trivially a bounded Donsker class, $\{\beta'Z(t), \beta \in \overline{B}_0, t \in [0, \tau]\}$ is Donsker. Since $\exp(\cdot)$ is Lipschitz on compacts and $\{Y(t), t \in [0, \tau]\}$ is monotone and bounded, the class $\{Y(t)e^{\beta'Z(t)}, \beta \in \overline{B}_0, t \in [0, \tau]\}$ is Donsker. Because $\mathcal{A}^\varepsilon \equiv \{A : (\gamma, \beta, A) \in \Psi_\varepsilon\}$ is uniformly bounded in total variation, the map $Y(t)e^{\beta'Z(t)} \mapsto \int_0^\tau Y(s)e^{\beta'Z(s)} dA(s)$, as a map from an element in $\ell^\infty(\overline{B}_0 \times [0, \tau])$ to an element in $\ell^\infty(\overline{B}_0 \times \mathcal{A}^\varepsilon)$, is continuous and linear, and the continuous mapping theorem yields that $\int_0^\tau Y(s)e^{\beta'Z(s)} dA(s)$ is Donsker as a process in $\ell^\infty(\Psi_\varepsilon)$. By conditions (C) and (D1), $\dot{G}_\gamma(t)$, $\ddot{G}_\gamma(t)$, $G_\gamma^{(1)}(t)$, $\dot{G}_\gamma^{(1)}(t)$ and $[\dot{G}_\gamma(t)]^{-1}$ are Lipschitz in $\gamma$ and $t$ over the appropriate range. Thus,

$$\left[\delta \frac{\dot{G}_\gamma^{(1)}(H^\psi(V))}{\dot{G}_\gamma(H^\psi(V))} - G_\gamma^{(1)}(H^\psi(V))\right] \quad \text{and} \quad \left[\delta \frac{\ddot{G}_\gamma(H^\psi(V))}{\dot{G}_\gamma(H^\psi(V))} - \dot{G}_\gamma(H^\psi(V))\right]$$

are also Donsker as processes in $\ell^\infty(\Psi_\varepsilon)$. Similar results and the fact that both sums of Donsker classes and products of bounded Donsker classes are Donsker give the result. □

PROOF OF COROLLARY 1. We first prove (ii). Using arguments from the proof of Theorem 4 and applying the $Z$-estimator master theorem (Theorem 3.3.1 of VW) gives that $\sqrt{n}(\hat{\psi}_n^\circ - \psi_*) = \sqrt{n}\mathbb{P}_n^\circ U^\tau(\psi_*)(\sigma^{-1}(\cdot)) + o_P(1)$ unconditionally, where $o_P$ denotes a quantity approaching 0 in outer probability. Since $\sqrt{n}(\hat{\psi}_n^\circ - \hat{\psi}_n) = \sqrt{n}(\mathbb{P}_n^\circ - \mathbb{P}_n)U^\tau(\psi_*)(\sigma^{-1}(\cdot)) + o_P(1)$ unconditionally, (ii) follows by the multiplier central limit theorem (Theorem 2.9.6 of VW) since $U^\tau(\psi_*)(\sigma_{\psi_*}^{-1}(\cdot))$ is $P_0$-Donsker and, over this Donsker class,

$$\sqrt{n}(\mathbb{P}_n^\circ - \mathbb{P}_n) = \bar{\zeta}_n^{-1}\sqrt{n}\sum_{i=1}^n (\zeta_i - \bar{\zeta}_n)(\Delta_{X_i} - P_0)$$

$$= \sqrt{n}\sum_{i=1}^n (\zeta_i - 1)(\Delta_{X_i} - P_0) + o_p(1),$$

where $\bar{\zeta}_n = n^{-1}\sum_{i=1}^n \zeta_i$ and $\Delta_X$ is the point mass at $X$. Similar arguments establish (i), but the nonparametric bootstrap central limit theorem (Theorem 3.6.1 of VW) is used in place of Theorem 2.9.6. □

PROOF OF PROPOSITION 4. Since $P_0 \log(p_\psi/p_0) \leq - \int_\mathcal{X} (p_\psi^{1/2} - p_0^{1/2})^2 \times p_0 \, d\nu$ and $p_0 = p_{\psi_0}$, we are done if we can show that, for any $\psi \in \Psi$,

$$G_\gamma(H^\psi(t)) = G_{\gamma_0}(H^{\psi_0}(t)) \tag{9.11}$$

for all $t \in [0, \tau]$ implies $\psi = \psi_0$ almost surely. Note that this requirement is valid even for $\gamma \leq 0$, since by (C) and (D1) an appropriate extension of $G_\gamma$ and its corresponding density $p_\psi$ exist. Taking the derivative of both sides of (9.11) with respect to $t$ yields

$$\dot{G}_\gamma(H^\psi(t))e^{\beta' Z(t)} b(t) = \dot{G}_{\gamma_0}(H^{\psi_0}(t))e^{\beta_0' Z(t)}, \tag{9.12}$$

where $b \equiv a/a_0$. Letting $t \downarrow 0$ in (9.12) gives $b(0+)e^{\beta' Z(0+)} = e^{\beta_0' Z(0+)}$ by condition (D1). This implies $\beta = \beta_0$ since $\mathrm{var}[Z(0+)]$ is positive definite. Hence, $b(0+) = 1$ by condition (H). Setting $\beta = \beta_0$ and dividing both sides of (9.12) by $e^{\beta_0' Z(t)}$, differentiating with respect to $t$ and letting $t \downarrow 0$ gives $\ddot{G}_\gamma(0+)e^{\beta_0' Z(0+)} + \dot{b}(0+) = \ddot{G}_{\gamma_0}(0+)e^{\beta_0' Z(0+)}$, where $\dot{b} = db/dA_0$. Since $\mathrm{var}[\beta_0' Z(0+)] > 0$, $\dot{b}(0+) = 0$. This now proves $\gamma = \gamma_0$ since $\ddot{G}_\gamma(0+)$ is monotone and bounded in $\gamma$. Now $A = A_0$ follows trivially. □

PROOF OF THEOREM 5. With any $h \in H_r$ such that $\sigma_{\psi_0}(h) = 0$, define the regular parametric one-dimensional submodel $\psi_{0t}(h) \equiv \psi_0 + t\{h_1, h_2, \int_0^{(\cdot)} h_3(s) \, dA_0(s)\}$. Note that $\sigma_{\psi_0}(h) = 0$ implies

$$P_0 \left\{ \left. \frac{\partial^2}{(\partial t)^2} L_n(\psi_{0t}) \right|_{t=0} \right\} = P_0 \{U^\tau(\psi_0)(h)\}^2 = 0,$$

where the score operator $U^\tau$ is given by (3.4). But this implies $P_0 \{U^\tau(\psi_0)(h)|\mathcal{G}(n, y, t)\}^2 = 0$, where the random set $\mathcal{G}(n, y, t) \equiv \{N, Y : N(s) = n(s), Y(s) = y(s), s \in [t, \tau]\}$, has nonzero probability. This then implies that $U^t(\psi_0)(h) = 0$ almost surely for all $t \in [0, \tau]$ (here is where we need the dependence on $\tau$ mentioned above in Section 3.3). Assuming that $[\{N(s), Y(s), Z(s)\}, s \geq 0]$ is censored at $V \in (0, \tau]$,

$$\begin{aligned}
0 &= G_{\gamma_0}^{(1)}(H^{\psi_0}(t))h_1 \\
&\quad + \dot{G}_{\gamma_0}(H^{\psi_0}(t)) \int_0^t Y(s) e^{\beta_0' Z(s)} \{h_2' Z(s) + h_3(s)\} \, dA_0(s).
\end{aligned} \tag{9.13}$$

Taking the derivative with respect to $t$ and letting $t \downarrow 0$ yields $h_2' Z(0+) + h_3(0+) = 0$ since $\dot{G}_{\gamma_0}^{(1)}(0+) = 0$ and $\dot{G}_{\gamma_0}(0+) = 1$ by assumption (D1). But this implies $h_2 = 0$ by condition (A3). Dividing (9.13) by $\dot{G}_{\gamma_0}(H^{\psi_0}(t))$, differentiating with respect to $t$ and taking $h_2 = 0$ yields

$$\begin{aligned}
0 &= [\dot{G}_{\gamma_0}^{(1)}(H^{\psi_0}(t))\dot{G}_{\gamma_0}(H^{\psi_0}(t)) - G_{\gamma_0}^{(1)}(H^{\psi_0}(t))\ddot{G}_{\gamma_0}(H^{\psi_0}(t))]h_1 \\
&\quad + [\dot{G}_{\gamma_0}(H^{\psi_0}(t))]^2 h_3(t).
\end{aligned}$$

Differentiating again with respect to $t$ and letting $t \downarrow 0$ gives $0 = \ddot{G}_{\gamma_0}^{(1)}(0+) \times e^{\beta_0 Z(0+)} h_1 + \dot{h}_3(0+)$, where $\dot{h}_3 \equiv a_0^{-1} dh_3/dt$. Now (H) and (D1) yield $h_1 = 0$. Thus, (9.13) implies $h_3(t) = 0$ for all $t \in [0, \tau]$, and the desired result follows. □

PROOF OF THEOREM 6. Define
$$\overline{\Psi}_M = \{\psi = (\gamma, \beta, A) : \gamma \in [-\varepsilon_0(K_0 A(\tau)), m_1],$$
$$\beta \in \overline{B}_0, A \in \mathcal{A}_{(\gamma)}, 1/M \leq a \leq M\}.$$

For $h = (h_1, h_2, h_3)$, with $h_1 \in \mathbb{R}$, $h_2 \in \mathbb{R}^d$ and $h_3 \in L_2([0,\tau])$, also define the metric $\|h\|_{\{2\}} \equiv |h_1| + \sqrt{h_2^T h_2} + (\int_0^\tau h_3^2(s)\,ds)^{1/2}$ and let the space of all such $h$ with $\|h\|_{\{2\}} < \infty$ be denoted $H_{\{2\}}$. Let $P_\psi g = \int_\mathcal{X} g p_\psi \, d\nu$ and denote $\tilde{\sigma}_\psi = P_\psi \hat{\sigma}_\psi$. Arguments in the proofs of Theorems 4 and 5 can be readily reworked to yield that $\tilde{\sigma}_\psi$ is one-to-one and continuously invertible as an operator in $H_{\{2\}}$. Thus, by the uniform compactness of $\overline{\Psi}_M$, by continuity and by Proposition 4 and Theorem 5, there exist constants $b_1, b_2 > 0$ and $k_1, k_2, k_3 < \infty$ such that, for all $\psi \in \overline{\Psi}_M$,
$$\psi^h(\tilde{\sigma}_\psi(h)) \geq b_1 \|h\|_{\{2\}}^2 \quad \text{and} \quad \psi^h(\hat{\sigma}_\psi(h)) \leq k_1 \|h\|_{\{2\}}^2$$
$\nu$-almost surely $\forall h \in H_{\{2\}}$,
$$\|\tilde{\sigma}_\psi(h)\|_v \geq b_2 \|h\|_v \quad \text{and} \quad \|\hat{\sigma}_\psi(h)\|_v \leq k_2 \|h\|_v,$$
$\nu$-almost surely $\forall h \in H_\infty$,
and
$$1/k_3 \leq |p_\psi(X)| \leq k_3, \qquad \nu\text{-almost surely.}$$

Choose $\varepsilon_1 = (2/3)(\{b_1/k_1\} \wedge \{b_2/k_2\})$ and, for any $f \in D(\nu)$ [where $D(\nu)$ is defined in Section 6.1] and $\psi \in \Psi$, denote $\sigma_\psi^{(f)} \equiv \int_\mathcal{X} \hat{\sigma}_\psi p_{(f)} \, d\nu$. Then, for all pairs $(f, \psi)$ with $f \in D(\nu)$, $\psi \in \overline{\Psi}_M$ and $\int_\mathcal{X} |f - f_\psi| \, d\nu \leq \varepsilon_1$, we have $\psi^h(\sigma_\psi^{(f)}(h)) \geq (2/3) b_2 \|h\|_{\{2\}}^2$ for all $h \in H_{\{2\}}$ and $\|\sigma_\psi^{(f)}(h)\|_v \geq (2/3) b_2 \|h\|_v$ for all $h \in H_\infty$.

For the chosen $Q$, let $\varepsilon = \{\varepsilon_1^2/[20 k_3 Q(Q+k_3)]\} \wedge \{\varepsilon_1/3\}$. Then, for all pairs $(f, \psi)$ with $f \in D(\nu)$, $\psi \in \Psi_M$ and $\int_\mathcal{X} |f - f_\psi| \, d\nu \leq \varepsilon$, we have a unique maximizer
$$\psi_1 = \arg\max_{\psi_* \in \Psi_M} \int_\mathcal{X} \log(p_{\psi_*}/p_{(f)}) p_{(f)} \, d\nu$$
such that $\|\sigma_{\psi_1}^{(f)}(h)\|_v \geq (2/3) b_2 \|h\|_v$ for all $h \in H_\infty$. To see this, note that, for $(f, \psi)$ and all $\psi_*$ with $\int_\mathcal{X} |p_{\psi_*} - p_\psi| \, d\nu \leq (2/3)\varepsilon_1$, we have $\int_\mathcal{X} |p_{(f)} - p_{\psi_*}| \, d\nu \leq \int_\mathcal{X} |f - f_\psi| \, d\nu + \int_\mathcal{X} |p_\psi - p_{\psi_*}| \, d\nu \leq \varepsilon_1$. Thus, $\psi_*^h(\sigma_{\psi_*}^{(f)}(h)) \geq \psi_*^h(\tilde{\sigma}_{\psi_*}(h)) -$

$k_1\|h\|_{\{2\}}^2 \int_{\mathcal{X}} |p_{(f)} - p_{\psi_*}|\, d\nu \geq (1/3)b_1\|h\|_{\{2\}}^2$ for all $h \in H_{\{2\}}$ and, arguing in a similar manner, $\|\sigma_{\psi_*}^{(f)}(h)\|_v \geq (1/3)b_2\|h\|_v$ for all $h \in H_\infty$. Thus, $\int_{\mathcal{X}} \log(p_{\psi_*}/p_{(f)})p_{(f)}\, d\nu$ is a convex function in $\psi_*$ with a continuously invertible second derivative, provided $\int_{\mathcal{X}} |p_{\psi_*} - p_\psi|\, d\nu \leq (2/3)\varepsilon_1$.

Furthermore, whenever $\psi_*$ satisfies $\int_{\mathcal{X}} |p_\psi - p_{\psi_*}|\, d\nu \leq \varepsilon/10$, we have $\int_{\mathcal{X}} \log(p_{\psi_*}/p_{(f)})p_{(f)}\, d\nu \geq -k_3 Q \int_{\mathcal{X}} |p_{\psi_*} - p_{(f)}|\, d\nu \geq -k_3 Q(11/10)\varepsilon$, but whenever $\psi_*$ satisfies $\int_{\mathcal{X}} |p_\psi - p_{\psi_*}|\, d\nu \geq (2/3)\varepsilon_1$, we have $\int_{\mathcal{X}} \log(p_{\psi_*}/p_{(f)}) \times p_{(f)}\, d\nu \leq -[2(Q+k_3)]^{-1}(\int_{\mathcal{X}} |p_{(f)} - p_{\psi_*}|\, d\nu)^2 \leq -\varepsilon_1^2/[18(Q+k_3)] < -k_3 Q(11/10)\varepsilon$. Hence, any Kullback–Leibler maximizer $\psi_{(f)}$ must satisfy $\int_{\mathcal{X}} |p_\psi - p_{\psi_{(f)}}|\, d\nu \leq (2/3)\varepsilon_1$. The desired existence and uniqueness of $\psi_{(f)}$ now follow, and $P_{(f)}\hat{\sigma}_{\psi_{(f)}}$ is one-to-one since $\|P_{(f)}\hat{\sigma}_{\psi_{(f)}}(h)\|_v \geq b\|h\|_v$ for all $h \in H_\infty$ and some $b > 0$ by arguments given above. □

PROOF OF THEOREM 7. Assume without loss of generality that $\mathrm{E}[Z_1] = 0$. Note that $A_*$ is uniformly bounded and equicontinuous by Lemma 1. Define $\tilde{K}_\gamma(e,t) \equiv \dot{G}_\gamma(t) - e\ddot{G}_\gamma(t)/\dot{G}_\gamma(t)$, $\dot{K}_\gamma(e,t) \equiv \partial \tilde{K}_\gamma(e,t)/(\partial t)$ and $\tilde{U}(x;\psi) \equiv z_1\{e - H^\psi(v)\tilde{K}_\gamma(e, H^\psi(v))\}$, where $x = (v,e,z) \in \mathcal{X}$, $z = (z_1, z_2)$ and $z_1$ is a possible value of $Z_1$. For each $\psi = (\gamma, \beta = (\beta_1, \beta_2), A) \in \Psi$, define $\psi^{(0)} = (\gamma, (0, \beta_2), A)$.

Assume first that $F_0^{z_1}(\cdot)$ is constant in $z_1$. It is clear that $\mathrm{E}[\tilde{U}(X;\psi_*^{(0)})] = 0$. If $\beta_{*1} > 0$, then conditions (A1′), (A3′) and (E1) imply

$$\mathrm{E}[Z_1 H^{\psi_*}(V)\tilde{K}_{\gamma_*}(\delta, H^{\psi_*}(V))] > \mathrm{E}[Z_1 H^{\psi_*^{(0)}}(V)\tilde{K}_\gamma(\delta, H^{\psi_*^{(0)}}(V))] = 0,$$

and thus $\mathrm{E}[\tilde{U}(X;\psi_*)] < 0$, but this is a contradiction. Similar arguments show that if $\beta_{*1} < 0$, $\mathrm{E}[\tilde{U}(X;\psi_*)] > 0$, also yielding a contradiction. Hence, $\beta_{*1} = 0$ and (i) follows.

Now assume $F_0^{z_1}$ is strictly increasing in $z_1$. If $\beta_{*1} \leq 0$, then by condition (E1)

$$\mathrm{E}[Z_1 H^{\psi_*}(V)\tilde{K}_{\gamma_*}(\delta, H^{\psi_*}(V))] \leq \mathrm{E}[Z_1 H^{\psi_*^{(0)}}(V)\tilde{K}_{\gamma_*}(\delta, H^{\psi_*^{(0)}}(V))],$$

but

$$\mathrm{E}[Z_1\{\delta - H^{\psi_*^{(0)}}(V)\tilde{K}_{\gamma_*}(\delta, H^{\psi_*^{(0)}}(V))\}]$$

(9.14)
$$= \mathrm{E}[-Z_1 H^{\psi_*^{(0)}}(V)\dot{G}_{\gamma_*}(H^{\psi_*^{(0)}}(V))]$$

$$+ \mathrm{E}\left[Z_1\delta\left\{1 + \frac{H^{\psi_*^{(0)}}(V)\ddot{G}_{\gamma_*}(H^{\psi_*^{(0)}}(V))}{\dot{G}_{\gamma_*}(H^{\psi_*^{(0)}}(V))}\right\}\right].$$

Now both $\delta$ and $-V$ are stochastically increasing in $Z_1$. By condition (E1), $-t\dot{G}_{\gamma_*}(t)$ is strictly decreasing in $t$ and $1 + t\ddot{G}_{\gamma_*}(t)/\dot{G}_{\gamma_*}(t)$ is nonincreasing in $t$; thus (9.14) is positive. Hence $\mathrm{E}[\tilde{U}(X;\psi_*)] > 0$. This implies that $\beta_{*1} > 0$ and (ii) follows. A similar proof can be used to establish (iii). $\square$

PROOF OF PROPOSITION 5. If $\gamma_0 = 0$, the proof follows from standard results for the Cox model. Hence, assume $\gamma_0 > 0$. Without loss of generality also assume $\mathrm{E}[Z] = 0$. The score for $\beta$ is

$$(9.15) \quad \mathrm{E}\left[\int_0^\tau \left\{Z - \frac{\mathrm{E}[ZY(t)e^{\beta'Z}]}{\mathrm{E}[Y(t)e^{\beta Z}]}\right\} Y(t) e^{\beta_0'Z} \dot{G}_{\gamma_0}^\circ(H^{\psi_0}(t)) \, dA_0(t)\right].$$

Note that the derivative of this with respect to $\beta$ is negative definite; thus, a 0 of (9.15) would be the unique maximizer of the profile likelihood (profiling over $A$). Note that

$$\frac{\mathrm{E}[\beta_0' Z Y(t)]}{\mathrm{E}[Y(t)]} < \frac{\mathrm{E}[\beta_0' ZY(t)e^{\beta_0'Z}\dot{G}_{\gamma_0}^\circ(H^{\psi_0}(t))]}{\mathrm{E}[Y(t)e^{\beta_0'Z}\dot{G}_{\gamma_0}^\circ(H^{\psi_0}(t))]} < \frac{\mathrm{E}[\beta_0' ZY(t)e^{\beta_0'Z}]}{\mathrm{E}[Y(t)e^{\beta_0'Z}]},$$

since $\dot{G}_{\gamma_0}^\circ(H^{\psi_0}(t))$ is decreasing in $\beta_0'Z$ but $e^{\beta_0'Z}\dot{G}_{\gamma_0}^\circ(H^{\psi_0}(t))$ is increasing in $\beta_0' Z$ by condition (E1). Thus,

$$\mathrm{E}\left[\int_0^\tau \beta_*' \left\{Z - \frac{\mathrm{E}[ZY(t)e^{\beta_*'Z}]}{\mathrm{E}[Y(t)e^{\beta_*'Z}]}\right\} Y(t) e^{\beta_0'Z} \dot{G}_{\gamma_0}^\circ(H^{\psi_0}(t)) \, dA_0(t)\right] = 0$$

for $\beta_* = \alpha_1 \beta_0$, where $\alpha_1 \in (0, 1)$.

To evaluate the score at $\beta = \beta_*$ in directions orthogonal to $\beta_0$, let $\beta_1 = [I - R\beta_0(\beta_0' R \beta_0)^{-1}\beta_0']u$, where $R = \mathrm{var}[Z]$ and $u \in \mathbb{R}^d$. Then $\mathrm{E}[\beta_1' Z \int_0^\tau Y(t) e^{\beta_0'Z} \times \dot{G}_{\gamma_0}^\circ(H^{\psi_0}(t)) \, dA_0(t) | \beta_0'Z] = 0$ since $\mathrm{E}[\beta_1'Z|\beta_0'Z] = 0$ by condition (A3″). Similarly, $\mathrm{E}[\beta_1' ZY(t) e^{\beta_*'Z}] = 0$. Since this is true for any $u \in \mathbb{R}^d$, $\beta_* = \alpha_1 \beta_0$ is indeed the unique maximizer of the profile likelihood. $\square$

PROOF OF PROPOSITION 6. When $\gamma_0 = 0$, the result follows since the Cox model is a valid submodel for any of the proportional hazards frailty models, and consistency has been established in Proposition 4. Assume $\gamma_0 > 0$ and let $\tilde{K}_\gamma$ and $\dot{K}_\gamma$ be as defined above in the proof of Theorem 7. The expected score for $\beta$, profiling over $A$, now has the form

$$(9.16) \quad \mathrm{E}\left[\int_0^\tau \left\{Z - \frac{\mathrm{E}[ZY(t)e^{\beta'Z}\tilde{K}_\gamma(\delta, H^\psi(V))]}{\mathrm{E}[Y(t)e^{\beta'Z}\tilde{K}_\gamma(\delta, H^\psi(V))]}\right\} \right.$$
$$\left. \times Y(t) e^{\beta_0'Z} \dot{G}_{\gamma_0}^\circ(H^{\psi_0}(t)) \, dA_0(t)\right],$$

where $A(t)$ solves

$$A(t) = \int_0^t \frac{\mathrm{E}[Y(s)e^{\beta_0'Z}\dot{G}_{\gamma_0}^\circ(H^{\psi_0}(s))]}{\mathrm{E}[Y(s)e^{\beta'Z}\tilde{K}_\gamma(\delta, H^\psi(V))]} \, dA_0(s)$$

and where $A(\tau) < \infty$ by Lemma 1. Let $\beta_{(\alpha,c)} = \alpha\beta_0 + c\beta_1$, where $\beta_1 = [I - R\beta_1(\beta_0'R\beta_0)^{-1}\beta_0']u$, with $R = \text{var}[Z]$ and $u \in \mathbb{R}^d$, as in Proposition 5. Denote $\psi_{(\alpha,c)} = (\gamma, \beta_{(\alpha,c)}, A)$. After multiplying by $\beta_1'$, the expected score (9.16) becomes

$$
\begin{aligned}
(9.17) \quad &-\int_0^\tau \frac{\mathrm{E}[\beta_1'ZY(t)e^{\beta_{(\alpha,c)}'Z}\tilde{K}_\gamma(\delta, H^{\psi_{(\alpha,c)}}(V))]}{\mathrm{E}[Y(t)e^{\beta_{(\alpha,c)}'Z}\tilde{K}_\gamma(\delta, H^{\psi_{(\alpha,c)}}(V))]} \\
&\quad \times Y(t)e^{\beta_0'Z}\dot{G}_{\gamma_0}^\circ(H^{\psi_0}(t))\,dA_1(t).
\end{aligned}
$$

We now evaluate $g(\alpha, c) = \mathrm{E}[\beta_1'ZY(t)e^{\beta_{(\alpha,c)}'Z}\tilde{K}_\gamma(\delta, H^{\psi_{(\alpha,c)}}(V))]$. Note that $g(\alpha, 0) = 0$ by previous arguments, and

$$
\frac{\partial g(\alpha, c)}{\partial c} = \mathrm{E}[(\beta_1'Z)^2 Y(t)e^{\beta_{(\alpha,c)}'Z}\{\tilde{K}_\gamma(\delta, H^{\psi_{(\alpha,c)}}(V)) + H^{\psi_{(\alpha,c)}}(V)\dot{K}_\gamma(\delta, H^{\psi_{(\alpha,c)}}(V))\}] > 0,
$$

since $\tilde{K}_\gamma(t) + t\dot{K}_\gamma(t) > 0$ for all $0 \le t < \infty$ by condition (E1). This now implies (9.17) $< 0$, which means that the expected score (9.16) is positive in the direction $-\beta_1$.

Hence, $\beta_* = \alpha\beta_0$ for some $\alpha \in \mathbb{R}$. Note that if $\alpha \le 0$, then

$$
\mathrm{E}[\beta_0'ZY(t)e^{\alpha\beta_0'Z}\tilde{K}_{\gamma_*}(\delta, H^{\psi_*(\alpha,0)}(V))] \le \mathrm{E}[\beta_0'ZY(t)\tilde{K}_{\gamma_*}(\delta, H^{\psi_*(0,0)}(V))]
$$

by condition (E1). However, $\mathrm{E}[\beta_0'Z\{\delta - A(V)\tilde{K}_{\gamma_*}(\delta, H^{\psi_*(0,0)}(V))\}] > 0$ by arguments used in the proof of Theorem 7. Thus, (9.16) is strictly positive if $\beta_* = \alpha\beta_0$ and $\alpha \le 0$. Hence $\alpha > 0$. $\square$

PROOF OF PROPOSITION 7. By condition (E2),

$$
\begin{aligned}
\lim_{\gamma \downarrow 0} \gamma^{-1}\dot{G}_\gamma(t) &= \lim_{\text{var}[W]\downarrow 0} \frac{\mathrm{E}[We^{-Wt} - e^{-Wt}]}{\mathrm{E}[e^{-Wt}\text{var}[W]]} \\
&= \lim_{\text{var}[W]\downarrow 0} \frac{\mathrm{E}[[W-1](1 - [W-1]t)]}{\text{var}[W]} = -t
\end{aligned}
$$

and, arguing similarly,

$$
\begin{aligned}
\lim_{\gamma \downarrow 0} \gamma^{-1}G_\gamma(t) &= \lim_{\text{var}[W]\downarrow 0} \frac{\mathrm{E}[-\log(e^{-Wt}/e^{-t})]}{\text{var}[W]} \\
&= \lim_{\text{var}[W]\downarrow 0} \frac{-\log(1 + \text{var}[W]t^2/2)}{\text{var}[W]} = \frac{-t^2}{2}.
\end{aligned}
$$

By Proposition 5, the score test for $H_0: \gamma_0 = 0$ thus has limiting expectation

$$
-\mathrm{E}\left[\delta e^{\alpha_1\beta_0'Z}A_*(V) - \frac{e^{2\alpha_1\beta_0'Z}A_*^2(V)}{2}\right],
$$

where

$$A_*(t) = \int_0^t \frac{\mathrm{E}[Y(s)e^{\beta_0' Z}\dot{G}_{\gamma_0}^\circ(H^{\psi_*}(s))]}{\mathrm{E}[Y(s)e^{\alpha_1 \beta_0' Z}]}\, dA_0(s).$$

Since $A_*^2(V)/2 = \int_0^\tau Y(s) A_*(s)\, dA_*(s)$, the score expectation becomes

$$\begin{aligned}(9.18)\qquad &\mathrm{E}\bigg[-\int_0^\tau \bigg\{ e^{\alpha_1 \beta_0' Z} - \frac{\mathrm{E}[Y(t)e^{2\alpha_1 \beta_0' Z}]}{\mathrm{E}[Y(t)e^{\alpha_1 \beta_0' Z}]}\bigg\}\\ &\qquad\qquad \times Y(t) e^{\beta_0' Z} \dot{G}_{\gamma_0}^\circ(H^{\psi_0}(t)) A_*(t)\, dA_0(t)\bigg].\end{aligned}$$

This clearly is equal to 0 when $\gamma_0 = 0$. Validity under contiguous alternatives follows from the regularity of the estimators under the correct model as established in Theorem 5. □

**Acknowledgments.** We thank Robert J. Gray of Harvard University for providing the non-Hodgkin's lymphoma data. We also thank the editors and referees for several helpful suggestions which led to significant improvements in organization and presentation.

M. R. Kosorok  
Departments of Statistics  
 and Biostatistics and Medical Informatics  
University of Wisconsin  
1210 West Dayton Street  
Madison, Wisconsin 53706  
USA  
e-mail: kosorok@biostat.wisc.edu  

B. L. Lee  
Department of Statistics  
 and Applied Probability  
National University of Singapore  
3 Science Drive 2  
Singapore 117543  
Republic of Singapore  
e-mail: staleebl@nus.edu.sg

J. P. Fine
Departments of Statistics
  and Biostatistics and Medical Informatics
University of Wisconsin
1210 West Dayton Street
Madison, Wisconsin 53706
USA
e-mail: fine@biostat.wisc.edu